\documentclass[11pt,a4paper]{amsart}
\usepackage[T1]{fontenc}
\usepackage{graphicx,etoolbox,needspace,mathrsfs}
\usepackage{lmodern,amsmath,amssymb,amsthm,mathtools,booktabs,array,longtable}
\usepackage[margin=25mm]{geometry}
\usepackage[colorlinks=true,linkcolor=blue,citecolor=blue,urlcolor=blue]{hyperref}

\newcommand{\verificationurl}{https://github.com/luanyongzhi/BrauerSemisimplicity}

\newcommand{\C}{\mathbb C}\newcommand{\Q}{\mathbb Q}\newcommand{\Z}{\mathbb Z}
\newcommand{\Br}{\operatorname{Br}}\newcommand{\Irr}{\operatorname{Irr}}
\newcommand{\Ind}{\operatorname{Ind}}

\newcommand{\rank}{\operatorname{rank}}

\newtheorem{mainresult}{Theorem}

\newtheorem{theorem}{Theorem}[section]
\newtheorem{lemma}[theorem]{Lemma}
\newtheorem{proposition}[theorem]{Proposition}
\newtheorem{corollary}[theorem]{Corollary}

\newtheoremstyle{uprightdefinition}
  {6pt}{6pt}{\normalfont\rmfamily}{}{\bfseries}{.}{.5em}{}
\theoremstyle{uprightdefinition}
\newtheorem{definition}[theorem]{Definition}
\AtBeginEnvironment{definition}{\Needspace{7\baselineskip}}
\theoremstyle{remark}\newtheorem{remark}[theorem]{Remark}

\title[Semisimplicity Criteria for Brauer Algebras]{Semisimplicity of Brauer algebras associated with finite Coxeter types}

\author[S. Liu]{Shoumin Liu}
\address[S. Liu]{School of Mathematics, Shandong University, Jinan 250100, China}
\email{s.liu@sdu.edu.cn}
\author[Y. Luan]{Yongzhi Luan}
\address[Y. Luan]{School of Mathematics, Shandong University, Jinan 250100, China}
\email{luanyongzhi@email.sdu.edu.cn}


\date{\today}

\begin{document}

\begin{abstract}
A finite-dimensional complex algebra is semisimple if its Jacobson radical
is zero. In this paper, we determine the complex parameters for which
Brauer algebras associated with finite Coxeter types are semisimple.
We also determine semisimplicity at zero for the specified polynomial
diagram forms and show that it can depend on the choice of polynomial form.
\end{abstract}
\maketitle

\vspace{-\baselineskip}
\tableofcontents

\section*{Introduction}
\label{sec:intro}
Richard Brauer introduced the ordinary Brauer algebra in his study of
tensor representations of orthogonal and symplectic groups in 1937
\cite{Brauer1937}. Its diagram model $B_n(\delta)$ has a basis consisting
of the perfect matchings of $n$ top and $n$ bottom vertices. Multiplication is given
by diagram concatenation, with each closed middle loop replaced by a factor $\delta$.
Thus $\dim B_n(\delta)=(2n-1)!!$, independently of $\delta$. At
suitable integral parameters, its action on tensor space describes
the centralizer algebra of the corresponding classical group; in general, the abstract diagram
algebra and its representation image need not have the same dimension.

Later, Cohen, Frenk, Gijsbers, and Wales introduced and studied the simply laced extensions of Brauer algebras in
\cite{CGW2007,CGWDiagrams2007,CFW2009,CFWerr2019}. By diagram folding, one obtains further
families of Brauer algebras, including type $\mathscr{B}_n$ inside type $\mathscr{D}_{n+1}$
\cite{CohenLiu2014}, type $\mathscr{C}_n$ inside type $\mathscr{A}_{2n-1}$
\cite{CLY2012,Bowman2012}, and the folded exceptional and dihedral types
\cite{LiuF42013,LiuG22015,Liu2016}. Even when two Coxeter groups are isomorphic, the corresponding Brauer algebra models need not coincide. In particular, the $\mathscr{B}_n$ and
$\mathscr{C}_n$ models studied here have different contraction parameters, and the
$\mathscr{G}_2$ folding of $\mathscr{D}_4$ is distinct from the $\mathscr{I}_2(6)$ folding of $\mathscr{A}_5$.

We determine the semisimple complex fibers of these algebras, with the conventions in
Definitions~\ref{def:semisimplicity} and~\ref{def:polynomial-form}.
For a cellular algebra with cell forms $G_\lambda(\delta)$, the criterion of
Graham and Lehrer \cite[Theorem~3.8]{GL1996} gives
\begin{equation}\label{eq:cellular-test}
 A(\delta)\text{ is semisimple}
 \quad\Longleftrightarrow\quad
 \det G_\lambda(\delta)\ne0\quad\text{for every cell label }\lambda.
\end{equation}
We must therefore determine when all the cell forms are nondegenerate.

For the ordinary Brauer algebra $B_n(\delta)$, Rui and Rui--Si gave a necessary and
sufficient condition for semisimplicity in \cite{Rui2005,RuiSi2006}. For the symmetric type
$\mathscr{C}$ diagram algebra, Dimpi, Chowdhury, and Thangavelu gave a sufficient
condition and discussed its necessity in \cite[Theorem~5.7 and
Remark~5.8]{TypeC2026}. The generic decompositions in
\cite{CFW2009} and \cite[Theorems~4.1--4.2]{LiuMorita2025}
are over a rational function field. They do not directly determine
semisimplicity at each complex specialization.
We give the parameter conditions together with the structural reductions
and finite calculations that establish them.

For $n\ge2$, we define
\begin{equation}\label{eq:Zn}
 Z(n)=\{z\in\mathbb Z:4-2n\le z\le n-2\}
 \setminus\{z\in\mathbb Z:4-2n<z\le3-n,\ z\text{ odd}\},
\end{equation}
and put $2Z(n)=\{2z:z\in Z(n)\}$. For type $\mathscr{C}$ we use
\begin{equation}\label{eq:Q}
 Q_+(\delta)=\frac{\delta^2+\delta-2}{2},\qquad
 Q_-(\delta)=\frac{\delta^2-\delta}{2}.
\end{equation}
We list the finite sets $S_6,S_7,S_8$ in
Section~\ref{subsec:singleton-spectra} and define $S_F$ in
Section~\ref{subsec:type-F4}. We also define the dihedral polynomials $P_{m,q}$,
$Q_{k,q}^{(1)}$ and the matrices $T_{k,q}$ in
Section~\ref{subsec:type-I2}. We summarize our main results for nonzero parameters in the following theorem.

\clearpage
\begin{mainresult}[Semisimplicity at nonzero parameters]
\label{thm:master-nonzero-semisimplicity}
Let $\delta\in\mathbb C^\times$. Each Brauer algebra listed in
Table~\ref{tab:master-nonzero-semisimplicity} is semisimple if and only
if the corresponding condition holds, with the algebra model and loop normalization
given in the indicated section.
\end{mainresult}

\begin{longtable}{@{}>{\raggedright\arraybackslash}p{.14\textwidth}
 >{\raggedright\arraybackslash}p{.59\textwidth}
 >{\raggedright\arraybackslash}p{.20\textwidth}@{}}
\caption{Necessary and sufficient conditions for $\delta\ne0$.}
\label{tab:master-nonzero-semisimplicity}\\
\toprule
Type & Condition & Reference\\\midrule
\endfirsthead
\toprule Type & Condition & Reference\\\midrule
\endhead
\bottomrule\endfoot
$\mathscr{A}_{n-1}$, $n\ge2$ & $\delta\notin Z(n)$.
 & Theorem~\ref{thm:typeA-rui-si}\\[4pt]
$\mathscr{B}_n$, $n\ge2$ & $\delta\notin Z(n+1)\cup2Z(n)$.
 & Theorem~\ref{thm:Bn-fiber}\\[4pt]
$\mathscr{C}_1$ & Every nonzero parameter.
 & Section~\ref{sec:definitions}\\[4pt]
$\mathscr{C}_n$, $n\ge2$ & $Q_+(\delta)\notin Z(n)$ and $Q_-(\delta)\notin Z(n)$.
 & Theorem~\ref{thm:C-semisimple}\\[4pt]
$\mathscr{D}_n$, $n\ge4$ & $\delta\notin Z(n)\cup2Z(n)$.
 & Theorem~\ref{thm:D-complete}\\[4pt]
$\mathscr{E}_{\kappa}$, $\kappa=6,7,8$ & $\delta\notin S_{\kappa}$; the sets have $17,26,40$ elements.
 & Theorem~\ref{thm:app-exceptional}\\[4pt]
$\mathscr{F}_4$ & $\delta\notin S_F$: eight integers and twelve noninteger values.
 & Theorem~\ref{thm:f4-certified}\\[4pt]
$\mathscr{G}_2$ & $\delta\notin\{-2,1\}$ for the $39$-dimensional triality fold.
 & Proposition~\ref{prop:G2-criterion}\\[4pt]
$\mathscr{H}_3$ & $\delta^2\notin\{-12,-4,-1,1,3,4\}$.
 & Theorem~\ref{thm:app-h3}\\[4pt]
$\mathscr{H}_4$ & $\delta^2\notin\{-44,-16,-14,-11,-8,-4,1,4,8,11,13,16\}$.
 & Theorem~\ref{thm:H4-complete-criterion}\\[4pt]
$\mathscr{I}_2(m)$, odd $m\ge3$ & $P_{m,q}(\delta)\ne0$ for every $q\mid m$.
 & Theorem~\ref{thm:odd-dihedral-spectrum}\\[4pt]
$\mathscr{I}_2(4)$ & $\delta\notin\{-2,1\}$ for the natural four-strand diagram fold.
 & Proposition~\ref{prop:I24-C2}\\[4pt]
$\mathscr{I}_2(2k)$, $k\ge3$ & For $k=2^sn$ with $n$ odd,
$Q^{(1)}_{k,q}(\delta)\ne0$ for all $q\mid n$ and
$\det T_{k,q}(\delta)\ne0$ for all $q\mid k$.
 & Corollary~\ref{cor:even-dihedral-structure}\\
\end{longtable}

We recall the presentations and cellular preliminaries in
Section~\ref{sec:preliminaries}. For the classical families in
Section~\ref{sec:classical-cases}, we combine the ordinary Brauer criterion
with diagram ideals, colored diagram quotients, and an explicit
factorization of each type $\mathscr{C}$ cell form. For the exceptional types in
Section~\ref{sec:exceptional-cases}, we use admissible-set modules,
reflection class sums, and exact finite matrix calculations. For the
dihedral diagram families, we apply Fourier diagonalization to obtain criteria in
terms of Ramanujan sums. In Section~\ref{sec:further-results},
we determine semisimplicity at zero for the specified polynomial diagram
forms. We also construct two polynomial forms of the triality-folded
$\mathscr{G}_2$ algebra with different answers at zero.

\section{Definitions and cellular preliminaries}
\label{sec:preliminaries}

\subsection{Parameters, semisimplicity, and polynomial forms}
We use script letters for Coxeter types, such as $\mathscr{E}_6$, and
ordinary letters for generators, such as $E_6$. The group symbols $A_5$
and $C_2$ denote the alternating and cyclic groups, respectively.
Put $R=\mathbb C[z]$ and $R_\times=R[z^{-1}]$.

\begin{definition}\label{def:semisimplicity}
Let $A$ be a finite-dimensional unital complex algebra. Its
\textbf{Jacobson radical} $\operatorname{rad}A$ is the intersection of
all maximal left ideals of $A$. We say that $A$ is \textbf{semisimple} if
$\operatorname{rad}A=0$.
\end{definition}

By the Artin--Wedderburn theorem, $A$ is semisimple if and only if $A$ is isomorphic to a
finite direct product of full matrix algebras over $\mathbb C$, or equivalently, every
finite-dimensional $A$-module is a direct sum of simple modules.
A Laurent-parameter algebra has no fiber at zero. To study that parameter, we first choose a polynomial form.

\begin{definition}\label{def:polynomial-form}
Let $L$ be an $R_\times$-algebra. In this paper, a \textbf{polynomial form} of $L$
is a unital $R$-subalgebra $\mathcal A\subseteq L$, free of finite rank
over $R$, such that the natural map
$\mathcal A\otimes_RR_\times\longrightarrow L$ is an isomorphism.
For any $\delta\in\mathbb C$, we define the \textbf{fiber at $\delta$} by
$\mathcal A(\delta)=\mathcal A\otimes_{R,\,z\mapsto\delta}\mathbb C$.
For $L$ itself, we define its complex fiber by scalar extension from $R_\times$
only when $\delta\ne0$.
\end{definition}

Two polynomial forms of the same $R_\times$-algebra $L$ can have
nonisomorphic fibers at zero; we give an example in
Section~\ref{sec:further-results}. For nonzero parameters $\delta\in\mathbb C^\times$, we use the Laurent presentations.

\subsection{The Weyl-type presentations and their parameters}
We give the Weyl-type presentations over $R_\times$ and write $\delta$
for the specialization of $z$ in a complex fiber. The polynomial diagram
models will be introduced separately.
Each presentation below defines an associative unital algebra. One can also
define these algebras over any commutative ring with a chosen invertible
parameter by replacing $z$ with that parameter; throughout this paper, we consider their complex fibers.

\begin{definition}
\label{def:weyl-brauer}
\cite[Section~1, Table~1]{CFW2009},
\cite[Definition~5.1, equations~(21)--(44)]{LiuG22015}
For a connected finite Weyl diagram $M$, the \textbf{Brauer algebra}
$\Br(M;R_\times)$ is generated by $r_i$ and $e_i$ indexed by the nodes of $M$.
We use the node labels and exponents specified below.
At each node $i$ and for each pair of distinct nonadjacent nodes $i,j$,
we impose
\begin{align*}
 r_i^2&=1,&r_ie_i=e_ir_i&=e_i,&e_i^2&=z^{\nu_i}e_i,\\
 r_ir_j&=r_jr_i,&e_ir_j&=r_je_i,&e_ie_j&=e_je_i.
\end{align*}
For each ordered pair $(i,j)$ joined by a single edge, we impose
\[
 r_ir_jr_i=r_jr_ir_j,\qquad
 r_jr_ie_j=e_ie_j,\qquad r_ie_jr_i=r_je_ir_j.
\]
We set the exponents $\nu_i$ according to the following table:
\[
\begin{array}{c|l}
 M&\nu_i\\\hline
 \mathscr{A},\mathscr{D},\mathscr{E}&1\text{ at every node}\\
 \mathscr{B}_n&\nu_0=2,\quad\nu_i=1\ (1\le i<n)\\
 \mathscr{C}_n&\nu_0=1,\quad\nu_i=2\ (1\le i<n)\\
 \mathscr{F}_4&\nu_1=\nu_2=2,\quad\nu_3=\nu_4=1\\
 \mathscr{G}_2&\nu_1=1,\quad\nu_2=3.
\end{array}
\]
For a double edge, we take the ordered pair
$(i,j)=(1,0)$ in $\mathscr{B}_n$, $(0,1)$ in $\mathscr{C}_n$, and $(3,2)$ in $\mathscr{F}_4$, and we impose
the additional relations
\begin{align*}
 r_jr_ir_jr_i&=r_ir_jr_ir_j,& r_jr_ie_j&=r_ie_j,\\
 r_je_ir_je_i&=e_ie_je_i,& (r_jr_ir_j)e_i&=e_i(r_jr_ir_j),\\
 e_jr_ie_j&=ze_j,& e_je_ie_j&=ze_j,\\
 e_jr_ir_j&=e_jr_i,&e_je_ir_j&=e_je_i.
\end{align*}
For the triple edge of $\mathscr{G}_2$, we take $(i,j)=(2,1)$ and impose
\begin{align*}
 r_ie_je_i&=r_je_i,&e_ie_jr_i&=e_ir_j,\\
 e_jr_ie_jr_ie_j&=e_j,&e_jr_ie_jr_ir_j&=e_jr_ir_jr_i,\\
 e_ir_je_i&=z^2e_i,&r_jr_ie_jr_ie_j&=r_ir_jr_ie_j,\\
 (r_jr_i)^6&=1.
\end{align*}
Thus our $\mathscr{G}_2$ indices $(1,2)$ correspond to indices $(1,0)$ in
\cite{LiuG22015}. In $\mathscr{B}_n,\mathscr{C}_n$ the remaining nodes form the chain
$1,\ldots,n-1$; in $\mathscr{F}_4$ the chain is $1,2,3,4$.
For $\mathscr{C}_1$, the single node has exponent one.
\end{definition}

By \cite[Proposition~5.3]{LiuG22015}, the map that fixes all generators and reverses words gives an
$R_\times$-linear anti-involution, which we denote by $*$. If we set all $e_i=0$, the relations in Definition~\ref{def:weyl-brauer} reduce to
the Coxeter relations of $W(M)$, so the quotient by the two-sided ideal generated by the $e_i$
is $R_\times[W(M)]$. The map sending the simple reflections to the generators
$r_i$ splits this quotient map, and thus the Coxeter group algebra
embeds into $\Br(M;R_\times)$. After specializing at $\delta\ne0$, the quotient $\mathbb C[W(M)]$
is semisimple. We therefore reduce semisimplicity to the ideal generated
by the $e_i$.

In type $\mathscr{A}$, we also use the diagram model, which gives a natural
polynomial form over $\mathbb C[z]$ with a well-defined fiber at $z=0$.

\begin{definition}\label{def:ordinary-brauer}
\cite[Section~4, equation~(4.2)]{GL1996}
A \textbf{Brauer diagram} on $n$ strands is a perfect matching of
$n$ top and $n$ bottom vertices. We call an edge joining the top and bottom rows a
\textbf{through strand}, and an edge within a single row a \textbf{horizontal arc}.
We call the number of through strands the \textbf{through rank}. The \textbf{ordinary Brauer diagram
algebra} $B_n(z)$ is the free $\mathbb C[z]$-module spanned by these diagrams.
To multiply two diagrams, we stack the first diagram above the second, replace each
closed middle component by a factor $z$, and keep the induced
matching of the external vertices. We define the anti-involution $*$ by reflecting a diagram across a horizontal line.
We set $B_0(z)=\mathbb C[z]$ and $(-1)!!=1$, and we denote the complex fiber at $z=\delta$
by $B_n(\delta)$.
\end{definition}

The rank of $B_n(z)$ over $\mathbb C[z]$, as well as the dimension of every fiber $B_n(\delta)$, is $(2n-1)!!$.
Here the identity diagram has $n$ vertical strands, $r_i$ interchanges the
adjacent strands $i$ and $i+1$, and $e_i$ joins positions $i,i+1$ horizontally in
both rows while keeping all other strands vertical. These diagrams give
the realization of type $\mathscr{A}_{n-1}$.
When an algebra serves as the ambient algebra for a folding construction, we write capital letters $R_i,E_i$ for its generators.
We will define the symmetric diagram model of type $\mathscr{C}$ in
Definition~\ref{def:symmetric-brauer}, whose basis consists of diagrams that are individually
invariant under reflection.

We also need the walled Brauer algebra when studying colored quotients in Section~\ref{sec:classical-cases}.
\begin{definition}\label{def:walled-brauer}
\cite[page~171, Section~2]{CDDM2008}
For integers $r,s\ge0$, we place a vertical wall after the first $r$ vertices in each row
of a Brauer diagram on $r+s$ strands. The \textbf{walled Brauer algebra}
$B^{\rm wall}_{r,s}(d)$ has a complex basis consisting of those Brauer diagrams whose
through strands stay on the same side of the wall and whose horizontal
arcs cross the wall. We define its multiplication by diagram concatenation with loop weight $d$.
\end{definition}

By folding, we embed the Brauer algebras of types $\mathscr{B}_n$, $\mathscr{C}_n$, $\mathscr{F}_4$, and $\mathscr{G}_2$ into those of types $\mathscr{D}_{n+1}$,
$\mathscr{A}_{2n-1}$, $\mathscr{E}_6$, and $\mathscr{D}_4$, respectively.

\begin{definition}\label{def:folded-realization}
For one of the specified partitions of the simply laced ambient diagram
into subsets $P$ of pairwise nonadjacent nodes, we set
\[
 r_P=\prod_{a\in P}R_a,\qquad e_P=\prod_{a\in P}E_a.
\]
We call the subalgebra generated by these products the \textbf{folded realization}
associated with the partition.
\end{definition}

Since the factors in each $P$ commute, we have $r_P^2=1$, $r_Pe_P=e_Pr_P=e_P$, and
$e_P^2=z^{|P|}e_P$, where the part size $|P|$ accounts for the exponent $\nu_i$ in Definition~\ref{def:weyl-brauer}.
For the four partitions used here, these folded realizations are isomorphic
to the corresponding Weyl-type algebras in Definition~\ref{def:weyl-brauer}, by
\cite[Theorem~2.6]{CohenLiu2014}, \cite[Theorem~1.1]{CLY2012},
\cite[Theorem~1.2]{LiuF42013}, and \cite[Theorem~6.1]{LiuG22015}.
Here we take the subalgebra generated by the products $r_P,e_P$.
It need not be the full fixed-point subalgebra of the linear diagram automorphism.

For type $\mathscr{C}_n$, we label the nodes of the ambient $\mathscr{A}_{2n-1}$ chain by $1,\ldots,2n-1$.
The partition consists of $\{n\}$ and $\{n-j,n+j\}$ for $1\le j<n$, so we have
$r_0=R_n$, $e_0=E_n$, $r_j=R_{n-j}R_{n+j}$, and
$e_j=E_{n-j}E_{n+j}$. For type $\mathscr{B}_n$, we label the two fork nodes of
$\mathscr{D}_{n+1}$ by $1,2$ (both adjacent to node $3$), followed by the chain
$3,4,\ldots,n+1$. Then we have
\[
 r_0=R_1R_2,\quad e_0=E_1E_2,\qquad
 r_j=R_{j+2},\quad e_j=E_{j+2}\quad(1\le j<n).
\]
When $n=2$, the ambient diagram is the three-node path $\mathscr{A}_3$. We will give the explicit embeddings for $\mathscr{F}_4$ and $\mathscr{G}_2$ in
Sections~\ref{subsec:type-F4} and~\ref{subsec:G2-certified}.

For the noncrystallographic types $\mathscr{H}_3,\mathscr{H}_4$ and the dihedral types $\mathscr{I}_2(m)$ (which are not included in Definition~\ref{def:weyl-brauer}), we use the folded realizations inside $\mathscr{D}_6$, $\mathscr{E}_8$, and $\mathscr{A}_{m-1}$ from \cite[Theorems~1.1--1.2]{Liu2016} and
\cite[Theorem~1.6]{LiuI22014}; we will recall their presentations and diagram
bases in Section~\ref{sec:exceptional-cases}.

\subsection{Cell data, cell modules, and bilinear forms}
We recall the cell modules and bilinear forms of Graham and Lehrer
\cite{GL1996}, using the lower-ideal convention throughout.

\begin{definition}\label{def:cellular}
\cite[Definition~1.1]{GL1996}
A \textbf{cell datum}
for an algebra $A_R$
over a commutative ring $R$ consists of a finite partially ordered set
$\Lambda$, finite nonempty sets $M(\lambda)$, an $R$-basis
\[
 \{C^\lambda_{s,t}:\lambda\in\Lambda,\ s,t\in M(\lambda)\},
\]
and an $R$-linear anti-involution $*$, such that
$(C^\lambda_{s,t})^*=C^\lambda_{t,s}$ and
\begin{equation}\label{eq:cell-axiom-expanded}
 aC^\lambda_{s,t}\equiv
 \sum_{u\in M(\lambda)}r_a(u,s)C^\lambda_{u,t}
 \pmod{A_R^{<\lambda}}.
\end{equation}
Here $A_R^{<\lambda}$ is spanned by the basis elements with labels
$\mu<\lambda$, and the coefficients $r_a(u,s)$ are independent of $t$.
An algebra admitting such a datum is a \textbf{cellular algebra}.
\end{definition}

A cell datum determines cell modules and their canonical bilinear forms.

\begin{definition}\label{def:cell-module-form}
\cite[Definitions~2.1, 2.3 and~3.1, Proposition~2.4]{GL1996}
For each $\lambda\in\Lambda$, the \textbf{cell module} $\Delta_R(\lambda)$ is the free $R$-module with basis $\{c_s:s\in M(\lambda)\}$ and action
$ac_s=\sum_u r_a(u,s)c_u$. Its \textbf{cell form} is the symmetric bilinear form
$\phi_\lambda$ defined by
\begin{equation}\label{eq:cell-form-expanded}
 C^\lambda_{s,t}C^\lambda_{u,v}
 \equiv\phi_\lambda(c_t,c_u)C^\lambda_{s,v}
 \pmod{A_R^{<\lambda}}.
\end{equation}
We define its \textbf{Gram matrix} by
$G_\lambda=(\phi_\lambda(c_s,c_t))_{s,t}$. Over a field, we define its
\textbf{form radical} by
$\operatorname{rad}\phi_\lambda=\{x:\phi_\lambda(x,y)=0
\text{ for every }y\}$.
\end{definition}

By the cellular axioms, we have
$\phi_\lambda(ax,y)=\phi_\lambda(x,a^*y)$,
so $\operatorname{rad}\phi_\lambda$ is an $A$-submodule of $\Delta(\lambda)$. Over $\mathbb C$, if
$\phi_\lambda$ is not identically zero, then the quotient
$L(\lambda)=\Delta(\lambda)/\operatorname{rad}\phi_\lambda$
is simple, and the nonzero quotients $L(\lambda)$ form a complete set of pairwise nonisomorphic simple $A$-modules
\cite[Proposition~3.2 and Theorem~3.4]{GL1996}. In this case, $\operatorname{rad}\phi_\lambda$
is the unique maximal submodule of $\Delta(\lambda)$. When $\phi_\lambda=0$,
we have $\operatorname{rad}\phi_\lambda=\Delta(\lambda)$, which need not coincide
with the Jacobson radical of $\Delta(\lambda)$.

We recall the cellular semisimplicity criterion
\cite[Theorem~3.8]{GL1996}.
\begin{proposition}[The cellular semisimplicity test]
\label{prop:cellular-test-expanded}
Let $A$ be a finite-dimensional cellular complex algebra and let
$G_\lambda$ be the Gram matrix of $\phi_\lambda$. Then we have
\[
 \dim(A/\operatorname{rad}A)=\sum_{\lambda\in\Lambda}
               (\operatorname{rank}G_\lambda)^2,
 \qquad
 \dim A=\sum_{\lambda\in\Lambda}|M(\lambda)|^2.
\]
In particular, $A$ is semisimple if and only if $G_\lambda$
is nonsingular for every $\lambda\in\Lambda$.
\end{proposition}
\begin{proof}
By the classification of simple $A$-modules recalled above and the Artin--Wedderburn theorem applied to $A/\operatorname{rad}A$, we obtain the first equality.
The second equality counts the cellular basis elements. Since
$\operatorname{rank}G_\lambda\le |M(\lambda)|$ for each $\lambda\in\Lambda$, the two
dimensions are equal if and only if every $G_\lambda$ has full rank. This proves
\eqref{eq:cellular-test}.
\end{proof}

\begin{remark}\label{rem:simple-cell-caution}
Equivalently, $A$ is semisimple if and only if all cell modules $\Delta(\lambda)$ are simple and pairwise nonisomorphic.
Pairwise nonisomorphism is necessary: for example, the algebra
$\mathbb C[x]/(x^2)$, with cellular basis $\{x,1\}$ ordered by $\lambda_x<\lambda_1$,
has two isomorphic one-dimensional cell modules whose cell forms are zero and nonzero, respectively, and $\mathbb C[x]/(x^2)$ is not
semisimple.
\end{remark}

We recall the composition-factor order from
\cite[Proposition~3.6]{GL1996}, using the lower-ideal convention in
\eqref{eq:cell-axiom-expanded}.
\begin{proposition}
\label{prop:cellular-order}
Let $\mu\in\Lambda$ be a label with $\phi_\mu\ne0$. If $L(\mu)$ is a
composition factor of $\Delta(\lambda)$, then we have
$\lambda\le\mu$. Moreover, if $\phi_\lambda\ne0$, then every composition factor
of $\operatorname{rad}\phi_\lambda$ has label $\mu>\lambda$.
\end{proposition}
\begin{proof}
We write $A^{\le\nu}$ for the two-sided ideal spanned by the cellular basis elements with labels at most $\nu$.
If $\lambda\not\le\mu$, then the intersection
$A^{\le\mu}\cap A^{\le\lambda}$ lies in $A^{<\lambda}$.
Hence we have $A^{\le\mu}A^{\le\lambda}\subseteq A^{<\lambda}$,
so $A^{\le\mu}$ annihilates $\Delta(\lambda)$.
On the other hand, since $\phi_\mu\ne0$, $A^{\le\mu}$ acts nontrivially on $L(\mu)$, which is a contradiction. This proves the first assertion.
By \eqref{eq:cell-form-expanded}, we also see that
$A^{\le\lambda}$ annihilates $\operatorname{rad}\phi_\lambda$, whereas
$A^{\le\lambda}$ acts nontrivially on $L(\lambda)$ when $\phi_\lambda\ne0$.
Therefore, $L(\lambda)$ cannot appear as a composition factor of $\operatorname{rad}\phi_\lambda$, which gives $\mu>\lambda$.
\end{proof}

In the diagram filtrations used here, increasing the number of horizontal contractions
moves down the cellular partial order. For the exceptional types, we will prove the corresponding cardinality
bound directly from the normal-form ideal filtration in Lemma~\ref{lem:specialized-layer-detection}.

Since scalar extension preserves a cell datum and its multiplication rules,
we can specialize a cellular basis over $\mathbb C[z]$ at any
$\delta\in\mathbb C$, whereas we can specialize a cellular basis over $R_\times$
only at $\delta\in\mathbb C^\times$. Under a change of basis with transition matrix $P$, a Gram matrix $G$ transforms into
$P^{\mathsf T}GP$, which has the same rank as $G$ whenever $P$
is invertible at the specialized parameter.

\subsection{Layer forms and elementary semisimplicity reductions}
When a normal-form basis gives an ideal filtration directly, we can use the following lemma to test semisimplicity on the successive layers.

\begin{lemma}\label{lem:sandwich-layers}
Suppose that a finite-dimensional complex algebra $A$ has a chain of two-sided ideals
$0=I_0\subset I_1\subset\cdots\subset I_\ell=A$
such that each successive quotient $I_r/I_{r-1}$ is isomorphic as a vector space to $M_m(H)$ for a finite-dimensional
semisimple coefficient algebra $H$ and an integer $m\ge1$, with multiplication given by
\[
 X\circ Y=X\Phi Y,\qquad \Phi\in M_m(H).
\]
Then $A$ is semisimple if and only if every coefficient block
$\Phi_\theta\in M_{m\dim\theta}(\mathbb C)$ is invertible
for every layer and every simple $H$-module $\theta$.
\end{lemma}
\begin{proof}
By the Artin--Wedderburn decomposition of $H$, each layer decomposes into a direct sum
of complex sandwich matrix blocks. If a block $\Phi_\theta$ is invertible, then the map $X\mapsto\Phi_\theta X$ is an algebra isomorphism onto the full
matrix algebra $M_{m\dim\theta}(\mathbb C)$.
If $\Phi_\theta$ is singular, then the block has no identity element, because any left identity $E$ would
satisfy $E\Phi_\theta=I$; hence it cannot be a finite-dimensional
semisimple algebra.

Since every two-sided ideal and every quotient of a finite-dimensional semisimple algebra is a direct product of full matrix algebras, the invertibility of each $\Phi_\theta$ is necessary. Conversely, any semisimple two-sided ideal $I$ of $A$ has a unit element $e$, which satisfies $ae=eae=ea$ for all $a\in A$. Thus $e$ is a central idempotent in $A$, so we have
$A=I\times(1-e)A$ with $(1-e)A\cong A/I$. By induction on the ideal
chain, we get the sufficiency.
\end{proof}

We distinguish the operator matrix $\Phi_\theta$ from the scalar Gram matrix $G_\theta$.

\begin{definition}\label{def:coefficient-forms}
Let $\mathcal I$ be the finite set of boundary labels in a given layer,
and let $\Phi=(\phi_{b,c})_{b,c\in\mathcal I}$ be its pairing matrix over
$\mathbb C[H]$. For a representation $\rho_\theta: \mathbb C[H]\to\operatorname{End}(U_\theta)$, we define the \textbf{normalized coefficient block} by
\[
 \Phi_\theta=(\rho_\theta(\phi_{b,c}))_{b,c\in\mathcal I}.
\]
Suppose that $U_\theta$ admits a nondegenerate symmetric $H$-invariant bilinear form
with matrix $S_\theta$, and that the layer anti-involution satisfies
$\phi_{b,c}^*=\phi_{c,b}$. With respect to the compatible cellular basis, we define the
\textbf{scalar layer form} by
\[
 \langle b\otimes u,c\otimes v\rangle
   =u^{\mathsf T}S_\theta\rho_\theta(\phi_{b,c})v,
 \qquad
 G_\theta=(I_{|\mathcal I|}\otimes S_\theta)\Phi_\theta.
\]
\end{definition}

From Definition~\ref{def:coefficient-forms}, we see that $G_\theta$ and $\Phi_\theta$ have the same rank, and their
determinants satisfy $\det G_\theta=(\det S_\theta)^{|\mathcal I|}\det\Phi_\theta$.
In a nonorthogonal basis of $U_\theta$, the matrix $\Phi_\theta$ need not be symmetric.
All coefficient groups $H$ and representations $U_\theta$ considered below carry such nondegenerate invariant forms; one can also test
semisimplicity directly from $\Phi_\theta$ via
Lemma~\ref{lem:sandwich-layers}.

The boundary set $\mathcal I$ of a cellular layer need not be a single Weyl-group orbit.
For instance, in type $\mathscr{F}_4$, the pair boundary set is a union of two orbits of sizes $12$ and $18$, and the
maximal boundary set is a union of two orbits of sizes $36$ and $3$; in each case, the entire union forms a single cellular
layer together with the mixed pairings between the two orbits.

We normalize each idempotent according to the actual square of the corresponding root product, using the simply laced formula $\delta^{-|X|}e_X$ when $X$ is an orthogonal set of ambient simply laced roots.
To apply Lemma~\ref{lem:sandwich-layers}, we check the layer multiplication,
all irreducible representations $\theta\in\Irr(H)$, and the dimension identity $\sum_{\theta\in\Irr(H)}(m\dim\theta)^2=m^2|H|$ across the full ideal chain.

The ideal splitting argument in Lemma~\ref{lem:sandwich-layers} also shows that
a finite-dimensional algebra $A$ is semisimple if and only if a two-sided ideal $I$
and the quotient $A/I$ are both semisimple. We will also use the corresponding
criterion for fiber products. If $A_+\to S$ and $A_-\to S$ are surjective
homomorphisms onto a semisimple algebra $S$, then $A_+\times_S A_-$ is
semisimple if and only if both $A_+$ and $A_-$ are semisimple.
Necessity follows from the two surjective projections onto $A_\pm$.
For sufficiency, the kernel ideals split off centrally, and the fiber product
is isomorphic to the direct product of those two kernels with $S$.

Finally, for any idempotent
$f\in A$, the truncation functor $M\mapsto fM$ is exact, so $fM=0$ for a finite-length $A$-module $M$ if and only if $fL=0$ for
every composition factor $L$ of $M$. We will use this property to analyze composition factors of specialized layer modules in Section~\ref{sec:exceptional-cases}.

\section{Classical types}
\label{sec:classical-cases}

We begin with the four classical infinite families $\mathscr{A}_{n-1}$, $\mathscr{B}_n$, $\mathscr{C}_n$, and $\mathscr{D}_n$.
Here $\mathscr{A}_{n-1}$ has rank $n-1$ and uses diagrams on $n$ strands, while
$\mathscr{C}_n$ has rank $n$ and uses symmetric diagrams on $2n$ strands.
The folded $\mathscr{B}_n$ algebra, the symmetric $\mathscr{C}_n$ diagram algebra, and the
simply laced $\mathscr{D}_n$ algebra have different diagram structures and loop normalizations.
We therefore treat their reductions to the ordinary Brauer algebra
$B_n(\delta)$ separately.

\subsection{Type \texorpdfstring{$\mathscr{A}_{n-1}$}{A(n-1)}: Rui and Rui--Si's classical criterion}
\label{subsec:type-A}

In type $\mathscr{A}_{n-1}$ ($n\ge2$), the ordinary Brauer algebra $B_n(\delta)$ has
dimension $b_n=(2n-1)!!$ and satisfies $e_{\kappa}^2=\delta e_{\kappa}$ for $1\le \kappa\le n-1$.
By \cite[Definition~4.9 and Theorem~4.10]{GL1996}, its cellular structure
is indexed by pairs $(a,\lambda)$, where
$0\le a\le \lfloor n/2\rfloor$ is the number of horizontal arcs in each row and
$\lambda\vdash n-2a$ is a partition of the number of through strands.

\begin{definition}\label{def:ordinary-pairing}
An \textbf{$a$-arc half-diagram} is a partial matching of $n$ ordered
vertices with $a$ horizontal arcs and $n-2a$ free vertices, where the free vertices are listed in
their natural order. We denote by $D_{n,a,\lambda}(z)$ the
\textbf{normalized ordinary pairing matrix}
obtained by overlaying two $a$-arc half-diagrams and evaluating the
resulting $\mathbb C[z][\mathfrak S_{n-2a}]$-valued pairing on the Specht module $S^\lambda$.
Here we transport the ordered free vertices from the column half-diagram to the row half-diagram; any overlay that decreases the
through rank evaluates to zero, and each closed middle loop contributes a factor $z$.
\end{definition}

If $H_\lambda$ is the Gram matrix of a fixed nondegenerate $\mathfrak S_{n-2a}$-invariant bilinear form on $S^\lambda$,
then the scalar cellular Gram matrix is $(I\otimes H_\lambda)D_{n,a,\lambda}$.
By Proposition~\ref{prop:cellular-test-expanded}, $B_n(\delta)$ is semisimple if and only if
every $D_{n,a,\lambda}(\delta)$ is nonsingular. We recall the classical semisimplicity criterion of Rui and Rui--Si as follows.

\begin{theorem}[{Rui \cite{Rui2005}; Rui--Si \cite[Theorem~2.3 and Lemma~2.4]{RuiSi2006}}]
\label{thm:typeA-rui-si}
Over $\mathbb C$, for $n\ge2$ and $\delta\ne0$, the ordinary Brauer algebra
$B_n(\delta)$ is semisimple if and only if $\delta\notin Z(n)$, where $Z(n)$ is
defined in \eqref{eq:Zn}.  At $\delta=0$, $B_n(0)$ is semisimple if and only if
$n\in\{1,3,5\}$ (with $B_1(\delta)=\mathbb C$ semisimple for all $\delta\in\mathbb C$).
\end{theorem}

We note that $0\in Z(n)$ for every $n\ge2$, and $Z(v)\subseteq Z(n)$ whenever
$2\le v\le n$. In our reductions for types $\mathscr{B}_n, \mathscr{C}_n$, and $\mathscr{D}_n$ at $\delta\ne0$,
the value $0\in Z(n)$ arises from even-strand ordinary factors
$B_2(0)$ evaluated at shifted parameters (such as $Q_\pm(\delta)=0$), so the zero-parameter exceptions $n\in\{3,5\}$ of $B_n(0)$ do not affect the nonzero-parameter criteria.

\subsection{Type \texorpdfstring{$\mathscr{B}_n$}{Bn}}
\label{subsec:type-B}
\label{subsec:Bn-reduction}
We determine semisimplicity for type $\mathscr{B}_n$ using an ordinary horizontal
ideal and two colored diagram quotients. Throughout this subsection, we assume $n\geq2$ and $\delta\in\mathbb C^\times$,
and we work with the Cohen--Liu algebra $\operatorname{Br}(\mathscr{B}_n;\delta)$ from \cite[Definition~2.1]{CohenLiu2014},
which satisfies $e_0^2=\delta^2e_0$ and $e_{\kappa}^2=\delta e_{\kappa}$ for $\kappa>0$.
Put $b_n=(2n-1)!!$. We will define the two-color diagram algebras
$C_\epsilon(n,d)$ and establish their basic properties independently in
Section~\ref{sec:D-reduction}.

\begin{theorem}\label{thm:Bn-fiber}
Let $n\ge2$, $\delta\in\mathbb C^\times$, $A=\operatorname{Br}(\mathscr{B}_n;\delta)$, and
$\Theta=Ae_0A$. Then we have algebra isomorphisms
\begin{align}
 \Theta&\cong I_{n+1}(\delta),\label{eq:Bn-ideal}\\
 A/\Theta&\cong C_+(n,\delta/2)
       \mathop{\times}_{\mathbb C[W(\mathscr{B}_n)]}C_-(n,\delta/2),
       \label{eq:Bn-fiber}
\end{align}
where $I_{n+1}(\delta)$ is the horizontal ideal of the ordinary Brauer
algebra $B_{n+1}(\delta)$. In particular,
\[
 A\text{ is semisimple}\quad\Longleftrightarrow\quad
 \delta\notin Z(n+1)\cup2Z(n).
\]
\end{theorem}
\begin{proof}
The diagram basis of $A$ given in
\cite[Notation~3.3 and Theorem~3.6]{CohenLiu2014} is
\[
 T_{n+1}^{|}\ \sqcup\ \xi T_{n+1}^{|=}
       \ \sqcup\ \theta U_{n+1}^{=}.
\]
In the first two families, the distinguished vertices $1$ and $\widehat1$ are
joined by a through strand. The third family consists of all undecorated
horizontal matchings multiplied by $\theta$. From the normal forms in
\cite[Notation~4.7, Table~1, and Theorem~5.7]{CohenLiu2014}, we see that the span of $\theta U_{n+1}^{=}$ is precisely $\Theta=Ae_0A$:
indeed, the normal-form families with indices $i=3,\ldots,6$ contain $e_0$
explicitly, and the mate-pair family with index $i=2$ lies in $Ae_0A$ by \cite[Lemma~4.5, equation~(4.4)]{CohenLiu2014}.
In fact, $\theta U_{n+1}^{=}$ spans the entire $\theta$-ideal
of the ambient $\mathscr{D}_{n+1}$ diagram algebra, with the same multiplication law.
Therefore, for $n\ge3$, we get \eqref{eq:Bn-ideal} from \eqref{eq:D-theta-ordinary}. When $n=2$, we also get $\Theta\cong I_3(\delta)$ by the same diagram multiplication, since
$\theta U_3^{=}$ consists of the nine horizontal diagrams of $B_3(\delta)$ and
$\delta^{-1}e_0$ corresponds to the ordinary contraction.

Next, for $\epsilon\in\{+1,-1\}$, we define the quotient maps
\[
 q_\epsilon:A\longrightarrow C_\epsilon(n,\delta/2),\qquad
 r_0\longmapsto T_1,\quad e_0\longmapsto0,\quad
 r_{\kappa}\longmapsto S_{\kappa},\quad e_{\kappa}\longmapsto E_{\kappa}\quad(\kappa>0).
\]
By diagram concatenation, the images of $r_\kappa,e_\kappa$ ($\kappa>0$) satisfy the ordinary Brauer relations. For the relations at the
double edge, we check directly that
\[
 T_1E_1T_1E_1=0,
 \qquad (T_1S_1T_1)E_1=\epsilon E_1
                   =E_1(T_1S_1T_1),
 \qquad (T_1S_1)^2=(S_1T_1)^2.
\]
All other defining relations involving $e_0$ are mapped to zero. Hence $q_+$ and $q_-$ are well-defined algebra
homomorphisms. Since $T_1$ and the permutations $S_\kappa$
generate all $T_{\kappa}$, and the word projectors together with the contractions $E_\kappa$ generate
the colored diagram basis of $C_\epsilon(n,\delta/2)$, we see that $q_\epsilon$ is surjective.

Let $J$ be the horizontal ideal of $A$ generated by $e_0$ and $e_1$. On $J/\Theta$, the central
bimodule operator $\xi$ satisfies $\xi^2=\delta^2$ and the mate-reflection relation
$\xi e_1=\delta(r_0r_1r_0)e_1$. Since $q_\epsilon(e_0)=0$, we obtain
$q_\epsilon(\xi z)=\epsilon\delta q_\epsilon(z)$ for all $z\in J$.
From the basis pairs $\{t,\xi t\}$, we see that the two $\xi$-eigenspaces in $J/\Theta$ have the same
dimension $h=2^n(b_n-n!)$. We also have
\[
 \dim(A/\Theta)=2^{n+1}b_n-2^n n!,\qquad
 \dim C_\epsilon(n,\delta/2)=2^n b_n.
\]
Therefore, the kernel of the induced map $\bar q_\epsilon: A/\Theta\to C_\epsilon(n,\delta/2)$ is precisely the
$-\epsilon\delta$-eigenspace of $\xi$ on $J/\Theta$, and the two kernels intersect trivially.
Modulo the horizontal ideals, both $q_+$ and $q_-$ induce the identity on the group algebra $\mathbb C[W(\mathscr{B}_n)]$. Comparing dimensions, we obtain \eqref{eq:Bn-fiber}.

By the ideal, quotient, and fiber-product reductions in
Section~\ref{sec:D-reduction}, we have
\[
 A\text{ is semisimple}\quad\Longleftrightarrow\quad
 B_{n+1}(\delta),\ C_+(n,\delta/2),\text{ and }C_-(n,\delta/2)
 \text{ are semisimple}.
\]
By Theorem~\ref{thm:typeA-rui-si} and the criterion for $C_+(n,\delta/2)$ in Section~\ref{sec:D-reduction}, the first two factors are semisimple if and only if
$\delta\notin Z(n+1)$ and $\delta\notin2Z(n)$, respectively.
For $n\ge4$, $C_-(n,\delta/2)$ is semisimple whenever
$\delta\notin2\{2-n,\ldots,n-2\}$, and we show in
Section~\ref{sec:D-reduction} that $2\{2-n,\ldots,n-2\}\subseteq Z(n+1)\cup2Z(n)$.
When $n=2$, the only nontrivial walled
factor in $C_-(2,\delta/2)$ is $B^{\rm wall}_{1,1}(\delta/2)$, which is semisimple for all
$\delta\ne0$. When $n=3$, the only nontrivial walled factors are
$B^{\rm wall}_{1,2}(d)$ and $B^{\rm wall}_{2,1}(d)$, where $d=\delta/2$.
Their horizontal cell pairing is
$\left(\begin{smallmatrix}d&1\\1&d\end{smallmatrix}\right)$,
so their nonzero exceptional parameters are $\delta=\pm2$.
Both lie in $Z(4)\cup2Z(3)$. This proves the criterion in every rank. In particular, the sets of nonzero nonsemisimple parameters for
$\mathscr{B}_2$ and $\mathscr{B}_3$ are $\{-2,1\}$ and
$\{-4,-2,1,2\}$, respectively.
\end{proof}

\subsection{Type \texorpdfstring{$\mathscr{C}_n$}{Cn}}
\label{subsec:type-C}
\label{sec:definitions}
\label{sec:gram}
We use the diagram basis for type $\mathscr{C}$ to factor every layer form
into ordinary Brauer forms, including at shifted loop parameters that may vanish.

\begin{definition}\label{def:symmetric-brauer}
\cite[Definition~2.5]{CLY2012}
On $2n$ top and $2n$ bottom vertices, we pair each vertex $\kappa$ with its mirror image $\bar\kappa$. The \textbf{symmetric Brauer diagram algebra}
$A_n(t)$ is the free $\mathbb C[t]$-module spanned by the Brauer diagrams that are
individually invariant under simultaneous reflection of the top and bottom rows.
We define its multiplication by ordinary diagram concatenation with loop weight $t$,
and we denote its complex fiber at $t=\delta$ by $A_n(\delta)$.
\end{definition}

Since diagram concatenation preserves reflection symmetry, $A_n(t)$ is
a $\mathbb C[t]$-subalgebra of $B_{2n}(t)$. For $n\ge2$, $A_n(t)$ is a proper subalgebra of the full
fixed-point subalgebra of the reflection automorphism on $B_{2n}(t)$, as the latter also
contains symmetric sums of non-invariant diagram pairs. For $\delta\ne0$, $A_n(\delta)$ is isomorphic to
the Cohen--Liu--Yu algebra $\Br(\mathscr{C}_n;\delta)$ by
\cite[Theorem~1.1]{CLY2012}, where
$e_0^2=\delta e_0$ and $e_\kappa^2=\delta^2e_\kappa$ for $\kappa>0$.

Using the mirror quotient, we represent each symmetric diagram by
signed edges and occupied singletons.
\begin{definition}\label{def:mirror-quotient}
The \textbf{mirror quotient} of a symmetric diagram is the induced partial
matching on its $2n$ mirror orbits of vertices. Under this quotient, a self-conjugate arc becomes an
\textbf{occupied singleton}, and a mirror pair of conjugate arcs becomes a
\textbf{signed edge}. Choosing the two sides of the mirror to label the two
vertices in each orbit, we assign the sign $+1$ to a signed edge if its two lifted
arcs connect vertices on the same side of the mirror, and $-1$
if they connect vertices on opposite sides.
\end{definition}

A closed signed cycle lifts
to two ordinary loops when its sign product is positive and one loop when
it is negative. A path between two occupied singleton ends lifts to one loop.
Counting the signed quotient matchings gives
\begin{equation}\label{eq:dimension}
 \dim A_n(\delta)=\sum_{j=0}^n\frac{(2n)!}{(2n-2j)!j!},
\end{equation}
where the $j$th summand counts the choice of $2j$ matched quotient vertices, their matching,
and their $2^j$ edge signs; in particular, $\dim A_n(\delta)$ is independent of the parameter.

In rank one,
\[
 A_1(\delta)=\C\langle s,F\rangle/
 (s^2-1,\ F^2-\delta F,\ sF-F,\ Fs-F)
\]
is semisimple exactly when $\delta\ne0$.

\subsubsection{Cell modules}
\begin{definition}\label{def:symmetric-half-diagrams}
\cite[Definitions~3.1--3.2]{Bowman2012}
A \textbf{symmetric $k$-half-diagram} is a reflection-invariant partial
matching on $2n$ vertices with $k$ horizontal arcs and $2(n-k)$ free vertices. In the mirror quotient, it has $s$ occupied
singletons, $j$ signed edges (with $s+2j=k$), and $n-k$ free
vertices. We denote the set of symmetric $k$-half-diagrams by $\mathcal D_{n,k}$ and keep the
sign convention of Definition~\ref{def:mirror-quotient}.
We write $W_m=(\mathbb Z/2\mathbb Z)\wr\mathfrak S_m$ for the
\textbf{hyperoctahedral group}. A \textbf{bipartition} of $m$
is an ordered pair of partitions $(\lambda^+,\lambda^-)$ with $|\lambda^+|+|\lambda^-|=m$.
\end{definition}

By direct counting, the cardinality of $\mathcal D_{n,k}$ is
\[
 m_{n,k}=\binom nk k!\sum_{j=0}^{\lfloor k/2\rfloor}
                     \frac1{(k-2j)!j!}.
\]
The cell labels of $A_n(\delta)$ are pairs $(k,\theta)$ with
$0\le k\le n$ and $\theta=(\lambda^+,\lambda^-)\vdash_2n-k$.
Put $p=|\lambda^+|$, $q=|\lambda^-|$, and
$d_\theta=\binom{p+q}{p}f^{\lambda^+}f^{\lambda^-}$, where
$f^\lambda=\dim S^\lambda$.
Then the cell module $\Delta_n(k,\theta)$ has dimension $m_{n,k}d_\theta$.

\Needspace{12\baselineskip}
\begin{definition}\label{def:type-c-pairing}
\cite[Lemma~3.3]{Bowman2012}
We fix a
matrix realization $\rho_\theta$ of the irreducible $W_{n-k}$-module
$U_\theta$ and extend it linearly to $\mathbb C[W_{n-k}]$. For two half-diagrams $b,c\in\mathcal D_{n,k}$, we overlay their arcs. We set
$\phi(b,c)=0$ if the number of free through strands decreases; otherwise, we set
$\phi(b,c)=\delta^\ell w$, where $\ell$ is the number of closed ordinary
loops in the lifted overlay and $w\in W_{n-k}$ maps the ordered free vertices
of the \emph{column} half-diagram $c$ to those of the \emph{row} half-diagram $b$ together with their signs.
We define the \textbf{normalized pairing matrix} by
\begin{equation}\label{eq:normalized-G}
 G_{n,k,\theta}(\delta)=\bigl(\rho_\theta(\phi(b,c))\bigr)_{b,c}.
\end{equation}
\end{definition}

The scalar cellular Gram matrix is obtained from \eqref{eq:normalized-G} by left multiplication
by the constant invertible Gram factor of the chosen nondegenerate invariant form on $U_\theta$.
Hence its rank agrees with that of \eqref{eq:normalized-G}, whereas its exact
determinant depends on the normalization of that form. Changing the matrix realization of $U_\theta$ uniformly
conjugates \eqref{eq:normalized-G} by similarity.

The diagram filtration by through strands has sections
\[
 \C\mathcal D_{n,k}\otimes\C W_{n-k}\otimes\C\mathcal D_{n,k},
\]
whose multiplication is given by the sandwich pairing $\phi$. Here diagram reflection interchanges
the two half-diagrams and inverts the group element. Inserting a cellular
basis of each complex group algebra, we obtain the cellular basis of $A_n(t)$
\cite[Lemmas~3.3 and~3.5]{Bowman2012},
\cite[Theorem~5.5 and (1.8)]{GL1996}. Since all multiplication coefficients are polynomials
in $t$ and no division by $t$ is used, this construction also defines
a cellular structure at $t=0$. From this cellular basis, we have the dimension identity
\[
 \dim A_n=\sum_{k=0}^n m_{n,k}^{2}\,2^{n-k}(n-k)!.
\]

By \eqref{eq:cellular-test}, semisimplicity is equivalent to
nondegeneracy of every matrix in \eqref{eq:normalized-G}.

The ordinary matrix $D_{v,a,\lambda}(z)$ of
Definition~\ref{def:ordinary-pairing} has order
\[
 d(v,a,\lambda)=\frac{v!}{(v-2a)!2^a a!}f^\lambda.
\]
In particular $D_{v,0,\lambda}=I_{f^\lambda}$ and
$D_{0,0,\varnothing}=(1)$. All normalizations in this section refer to these
matrices and to \eqref{eq:normalized-G}.

\subsubsection{Two explicit basis changes}
We factorize $G_{n,k,\theta}(\delta)$ in two steps: first we decouple the
occupied singletons, and then we diagonalize the action of the
edge-sign group.
For $\delta\ne0$, we apply the following substitution to each signed edge of a half-diagram:
\begin{equation}\label{eq:cut}
 e_\varepsilon(\kappa,j)\longmapsto
 e_\varepsilon(\kappa,j)-\delta^{-1}f_{\kappa}f_j,
\end{equation}
where $f_{\kappa}f_j$ denotes the pair of occupied singletons at $\kappa$ and $j$. This substitution preserves $k$, the active
vertex set, and the ordered set of free vertices. Expanding over all signed edges of a half-diagram, we obtain
a unitriangular transition matrix $P$ with respect to the number of signed edges, which
acts as the identity on the coefficient module $U_\theta$ and has determinant one.

\begin{definition}\label{def:signed-pairing}
We denote by $H_{m,j,\theta}(d_0,d_1)$ the
\textbf{signed diagram pairing} on signed $j$-edge half-diagrams
on $m$ quotient vertices with $m-2j=p+q$ free vertices, where a closed signed cycle
has weight $d_0$ if its sign product is $+1$ and $d_1$ if its sign product is $-1$, under the
same rank-drop and free-vertex transport rules as in
Definition~\ref{def:type-c-pairing}. This is a normalized coefficient
pairing: we evaluate the group-algebra entries in $\rho_\theta$, without
the additional coefficient Gram factor.
\end{definition}

The number of such signed $j$-edge half-diagrams is
$h_{m,j}=m!/((m-2j)!j!)$.
\begin{lemma}\label{lem:singleton-gram}
Under the basis change \eqref{eq:cut}, half-diagrams with distinct singleton sets are orthogonal.
For a fixed singleton set $S$ of size $s$ with $s+2j=k$, the corresponding diagonal block of
$P^{\mathsf T}GP$ is
\[
 \delta^s H_{n-s,j,\theta}(\delta^2-1,\delta-1).
\]
\end{lemma}
\begin{proof}
We examine each connected component in the overlay of two half-diagrams and expand the substitutions \eqref{eq:cut} by choosing which signed
edges to cut. Any path whose two free ends lie on the same
row reduces the through rank and hence contributes zero. For a path whose free ends lie on opposite rows,
any cut reduces the through rank, so only the uncut term survives and preserves the signed permutation of free vertices.
Similarly, any path joining a singleton to a free end reduces the through rank and contributes zero.
For a path of $u$ signed edges between two singleton ends, cutting $t$ edges produces
$t+1$ closed ordinary loops with coefficient $(-\delta^{-1})^t$. Summing over $t$, we get
$\delta\sum_{t=0}^u\binom ut(-1)^t=\delta(1-1)^u$.
Thus such a component survives if and only if $u=0$, which means the two half-diagrams share a common singleton at the same vertex, contributing a factor $\delta$.
Finally, a closed signed cycle has uncut weight $\delta^2$ (for positive sign product) or $\delta$ (for negative sign product).
Every nonempty choice of $t\ge1$ cuts on the cycle contributes $(-\delta^{-1})^t\delta^t=(-1)^t$,
and summing over $t\ge1$ gives $-1$. Hence the effective cycle weights become $\delta^2-1$ and $\delta-1$, respectively.
Multiplying the contributions over all components in the group algebra $\mathbb C[W_{n-k}]$ and applying $\rho_\theta$, we get our result.
\end{proof}

\begin{lemma}\label{lem:fourier-gram}
In the canonical induced realization of $U_\theta$, we choose the same Specht
bases as in the ordinary factors and apply the orthogonal edge basis change
$2^{-1/2}(e_++e_-)$ and $2^{-1/2}(e_+-e_-)$ to each signed edge.
Then $H_{m,j,\theta}(d_0,d_1)$ is congruent by a determinant-preserving transformation to the direct sum, over $a+b=j$, of
$\binom m{p+2a}$ copies of
\[
 2^j D_{p+2a,a,\lambda^+}\left(\frac{d_0+d_1}{2}\right)
       \otimes
       D_{q+2b,b,\lambda^-}\left(\frac{d_0-d_1}{2}\right).
\]
\end{lemma}
\begin{proof}
We assign the color $+1$ to both endpoints of $2^{-1/2}(e_++e_-)$ and the color $-1$ to both endpoints of $2^{-1/2}(e_+-e_-)$; the induced realization of $U_\theta$
assigns $+1$ to $p$ free vertices and $-1$ to $q$ free vertices.
If the colors disagree along an overlaid component, summing over the sign of the connecting edge gives zero.
For a monochromatic closed cycle with $u$ edges and color $\eta\in\{+1,-1\}$,
the character defined by $\chi_\eta(1)=1$ and
$\chi_\eta(-1)=\eta$ satisfies
\[
 \sum_{\varepsilon_1,\ldots,\varepsilon_u=\pm1}
 \chi_\eta(\varepsilon_1\cdots\varepsilon_u)
 d_{\varepsilon_1\cdots\varepsilon_u}
 =2^{u-1}(d_0+\eta d_1),
\]
where $d_{+1}=d_0$ and $d_{-1}=d_1$. Along a monochromatic through path of $u$ edges, the sign action at the free endpoints
cancels the character factors and leaves $2^u$ times the ordinary
permutation of same-color free vertices. Combining the $2j$ row and column edge sums with the normalization factor $2^{-j}$ from the basis change,
we get an overall factor $2^j$.
There are $\binom m{p+2a}$ choices of the $+1$-colored vertex set.
For each choice, the remaining basis is the tensor product of the ordinary
half-diagram and Specht bases for $\lambda^+$ and $\lambda^-$.
All Fourier coefficients are scalar and act as the identity on the
fixed Specht coordinates; no transposition of a coefficient representation
matrix is introduced. If we use an equivalent realization of $U_\theta$, we first conjugate uniformly by the corresponding similarity transformation. Since the orthogonal Fourier change on each half-diagram has
determinant $\pm1$, applying it on both rows and columns preserves the determinant.
\end{proof}

\subsubsection{All layers}
Combining the two basis changes above, we obtain the full decomposition of each cell form.
Here the multiplicities count the choices of singleton and color sets,
and the remaining tensor factors are ordinary Brauer pairing matrices.
For $\theta=(\lambda^+,\lambda^-)$ with $p+q=n-k$, we set
\[
 \mathcal I_k=\{(s,a,b)\in\Z_{\ge0}^3:s+2a+2b=k\},\qquad
 \nu_{s,a,b}=\binom ns\binom{n-s}{p+2a},
\]
and abbreviate
\[
 A_a^+=D_{p+2a,a,\lambda^+}(Q_+),\quad
 A_b^-=D_{q+2b,b,\lambda^-}(Q_-),\quad
 u_a=\dim A_a^+,\quad v_b=\dim A_b^-.
\]
\begin{theorem}[Factorization of the cell forms]\label{thm:master-gram}
For every $\delta\ne0$, the determinant-preserving invertible transformations in Lemmas~\ref{lem:singleton-gram} and~\ref{lem:fourier-gram} give the following block decomposition of $G_{n,k,\theta}(\delta)$:
\begin{equation}\label{eq:master-block}
 G_{n,k,\theta}(\delta)\ \longmapsto\
 \bigoplus_{(s,a,b)\in\mathcal I_k}
  [\delta^s2^{a+b}A_a^+\otimes A_b^-]^{\oplus\nu_{s,a,b}}.
\end{equation}
Taking determinants gives
\begin{equation}\label{eq:master-det}
 \det G_{n,k,\theta}
 =\prod_{(s,a,b)\in\mathcal I_k}
 \left[(\delta^s2^{a+b})^{u_av_b}
       (\det A_a^+)^{v_b}(\det A_b^-)^{u_a}\right]^{\nu_{s,a,b}}.
\end{equation}
Both formulas hold even when $Q_+=0$ or $Q_-=0$.
\end{theorem}
\begin{proof}
We first choose the singleton set $S$ of size $s$, then choose the $+1$-colored vertex set of size $p+2a$ on its complement,
and apply Lemmas~\ref{lem:singleton-gram} and~\ref{lem:fourier-gram}. The number of such choices is $\nu_{s,a,b}$, and the sum of the block sizes equals $m_{n,k}d_\theta$.
Applying the determinant formula for Kronecker products, we obtain \eqref{eq:master-det}.
Since our basis changes divide only by $\delta\ne0$ and by $2$, both \eqref{eq:master-block} and \eqref{eq:master-det} remain valid when $Q_+=0$ or $Q_-=0$.
\end{proof}

\subsubsection{The semisimplicity criterion}
\begin{theorem}\label{thm:C-semisimple}
For $n\ge2$ and $\delta\ne0$, we have
\[
 A_n(\delta)\text{ is semisimple}
 \quad\Longleftrightarrow\quad Q_+(\delta),Q_-(\delta)\notin Z(n).
\]
For $n=1$, $A_1(\delta)$ is semisimple if and only if $\delta\ne0$.
\end{theorem}
\begin{proof}
First, suppose that $Q_+(\delta)\in Z(n)\setminus\{0\}$.
By Theorem~\ref{thm:typeA-rui-si}, there exists a cell label $(a,\lambda)$ of $B_n(Q_+(\delta))$ such that $D_{n,a,\lambda}(Q_+(\delta))$ is singular.
Taking $k=2a$, $\theta=(\lambda,\varnothing)$, and $(s,a,b)=(0,a,0)$
in \eqref{eq:master-block}, we get a singular block in $G_{n,2a,\theta}(\delta)$.
If $Q_+(\delta)=0$, we take $k=n$, $\theta=(\varnothing,\varnothing)$, and
$(s,a,b)=(n-2,1,0)$, which gives the zero block
$2\delta^{n-2}D_{2,1,\varnothing}(0)=[0]$.
By symmetry, the same argument with $\theta=(\varnothing,\lambda)$ and $(s,a,b)=(0,0,b)$ or $(n-2,0,1)$ shows that $A_n(\delta)$ is not semisimple if $Q_-(\delta)\in Z(n)$.

Conversely, recall that $0\in Z(n)$ and $Z(v)\subseteq Z(n)$ for all $2\le v\le n$.
If $Q_+(\delta)\notin Z(n)$ and $Q_-(\delta)\notin Z(n)$, then $Q_\pm(\delta)\ne0$ and $Q_\pm(\delta)\notin Z(v)$ for every $2\le v\le n$. By Theorem~\ref{thm:typeA-rui-si} (together with the trivial cases $v=0,1$), all ordinary factors $D_{p+2a,a,\lambda^+}(Q_+)$ and $D_{q+2b,b,\lambda^-}(Q_-)$ in \eqref{eq:master-block} are invertible.
By \eqref{eq:cellular-test}, we conclude that $A_n(\delta)$ is semisimple.
\end{proof}

\begin{remark}\label{rem:typec2026-example59}
In \cite[Theorem~5.7]{TypeC2026}, the authors gave a sufficient condition for semisimplicity by excluding
integral values of both $\delta^2+\delta$ and $\delta^2-\delta$.
Theorem~\ref{thm:C-semisimple} gives a necessary and sufficient condition in finite rank: both
$Q_+(\delta)$ and $Q_-(\delta)$ must avoid $Z(n)$ (restricting
only one of these two shifted parameters is not sufficient).
For example, we have $Z(2)=\{0\}$, so $A_2(\delta)$ is semisimple at nonzero
parameters if and only if $\delta\notin\{-2,1\}$.
For $n=3$, $k=2$, and $\theta=((1),\varnothing)$, the cell form has order
$m_{3,2}=9$. By \eqref{eq:master-block}, it decomposes into a three-dimensional
singleton block, three negative-pair blocks, and one three-dimensional
positive-pair block, with determinant
\[
 \delta^6(2Q_-)^3(2Q_+-2)^2(2Q_++4)
 =\delta^9(\delta-1)^3(\delta^2+\delta-4)^2
       (\delta^2+\delta+2),
\]
which agrees with \cite[Example~5.9]{TypeC2026}. For nonzero parameters, this cell form is singular precisely when
$Q_-=0$ or $Q_+\in\{1,-2\}$; the full criterion uses all cells.
\end{remark}

\subsection{Type \texorpdfstring{$\mathscr{D}_n$}{Dn}}
\label{subsec:type-D}
\label{app:typeD}
\label{sec:D-reduction}

Throughout this subsection, we assume $n\geq4$ and $\delta\in\mathbb C^\times$, and we denote by
$A=\operatorname{Br}(\mathscr{D}_n;\delta)$ the Cohen--Frenk--Wales algebra of type $\mathscr{D}_n$.
We label the two fork nodes by $1,2$ (both adjacent to node $3$), followed by the chain
$3-4-\cdots-n$. Put $b_n=(2n-1)!!$ and recall that
\[
 Z(n)=\{4-2n,4-2n+1,\ldots,n-2\}
 \setminus\{j:4-2n\leq j\leq3-n,\ j\text{ odd}\}.
\]
Using two colored diagram quotients, we reduce the semisimplicity of
$\operatorname{Br}(\mathscr{D}_n;\delta)$ to that of ordinary and walled Brauer algebras.

\subsubsection{Two colored diagram algebras}

We describe the two quotients uniformly by using two rules for
the colors at a horizontal edge.

\begin{definition}\label{def:colored-algebras}
Put $d=\delta/2$ and let $\epsilon\in\{1,-1\}$.
The \textbf{two-color diagram algebra} $C_\epsilon(n,d)$ has as
objects the $2^n$ words of length
$n$ in the two colors $+1,-1$. A basis element from one word to another
is a Brauer matching of the upper and lower vertices, such that a through
strand joins equal colors and a horizontal strand joins colors whose
product is $\epsilon$. Composition is zero when the middle words differ.
When they agree, we concatenate the two diagrams, replace each closed colored component by a
factor $d$, and retain the external matching and colors. Thus $C_\epsilon(n,d)$
is the endomorphism algebra of the direct sum of the $2^n$ color words in the symmetric
diagram category where the dual of color $c$ is $\epsilon c$.
\end{definition}

Since the number of closed components
in a triple concatenation does not depend on the order of composition, this multiplication is associative.
Every uncolored matching has $2^n$ permitted colorings, with one free
choice on each of its $n$ strands. Therefore, we have
\[
 \dim C_\epsilon(n,d)=2^n b_n.
\]
Let $S_{\kappa}$ be the sum of the colored permutation diagrams for the adjacent
transposition $\kappa,\kappa+1$. Let $T_j$ act on a word by its color at position $j$.
Let $E_{\kappa}$ be the sum of all permitted colored cup--cap diagrams on positions
$\kappa,\kappa+1$, with vertical strands elsewhere. Their local relations include
\begin{gather*}
 T_j^2=1,\qquad S_{\kappa}T_{\kappa}S_{\kappa}=T_{\kappa+1},\qquad E_{\kappa}^2=2dE_{\kappa}=\delta E_{\kappa},
 \qquad S_{\kappa}E_{\kappa}=E_{\kappa}=E_{\kappa}S_{\kappa},\\
 T_{\kappa}E_{\kappa}=\epsilon T_{\kappa+1}E_{\kappa},\qquad
 E_{\kappa}T_{\kappa}=\epsilon E_{\kappa}T_{\kappa+1},\qquad E_{\kappa}T_{\kappa}E_{\kappa}=0.
\end{gather*}
Here the last equality holds because the two possible colors of the middle
cup contribute $d$ and $-d$. By diagram concatenation, the elements $S_{\kappa},E_{\kappa}$ satisfy the ordinary Brauer
relations with parameter $\delta$, and operators on disjoint
positions commute.
These elements generate $C_\epsilon(n,d)$: the products of
$(1\pm T_j)/2$ give the individual word idempotents, and multiplying
$S_{\kappa},E_{\kappa}$ by these idempotents yields the colored permutation and contraction
diagrams that generate every colored matching.

Let $\tau$ be the automorphism of $C_\epsilon(n,d)$ that reverses every color, so that
$\tau(T_j)=-T_j$, $\tau(S_{\kappa})=S_{\kappa}$, and $\tau(E_{\kappa})=E_{\kappa}$. We set
$F_\epsilon=C_\epsilon(n,d)^\tau$. Since $\tau$ fixes no
colored basis diagram, we have
\[
 \dim F_\epsilon=2^{n-1}b_n.
\]

\begin{lemma}\label{lem:D-colored-generators}
The algebra $F_\epsilon$ is generated by
\[
 S_1,\ldots,S_{n-1},\quad T_1T_2S_1,\quad
 E_1,\ldots,E_{n-1},\quad T_1E_1T_1.
\]
Moreover, $C_\epsilon(n,d)=F_\epsilon\rtimes\langle T_1\rangle$.
\end{lemma}
\begin{proof}
Let $B$ be the subalgebra generated by the displayed elements. We see that $B$ contains
$T_1T_2=(T_1T_2S_1)S_1$ and, by permutation conjugation, every $T_{\kappa}T_j$.
Moreover, conjugation by $T_1$ preserves $B$, since it interchanges $S_1$ with
$T_1T_2S_1$, interchanges $E_1$ with $T_1E_1T_1$, and fixes the other
listed generators and all $T_{\kappa}T_j$. If we adjoin $T_1$ to $B$, we get
every $T_j$ and hence the full algebra $C_\epsilon(n,d)$. Therefore, we have
$C_\epsilon(n,d)=B+BT_1$. Since every element of $B$ is fixed by $\tau$ and every
element of $BT_1$ lies in the $-1$-eigenspace of $\tau$, this sum is direct and
$B=F_\epsilon$. Since $T_1^2=1$, we obtain $C_\epsilon(n,d)=F_\epsilon\rtimes\langle T_1\rangle$.
\end{proof}

\begin{lemma}\label{lem:D-colored-map}
There is a surjective homomorphism $q_\epsilon:A\to F_\epsilon$ given by
\begin{align*}
 r_1&\longmapsto T_1S_1T_1=T_1T_2S_1,
 &e_1&\longmapsto T_1E_1T_1,\\
 r_{\kappa}&\longmapsto S_{\kappa-1},
 &e_{\kappa}&\longmapsto E_{\kappa-1}\qquad(2\leq \kappa\leq n).
\end{align*}
It kills $e_1e_2$ and satisfies $q_\epsilon(r_1e_2)=\epsilon E_1$.
\end{lemma}
\begin{proof}
The relations among vertices $2,\ldots,n$ are the ordinary Brauer
relations. The relations involving vertices $1,3$ are obtained from those
for $2,3$ by conjugation by $T_1$, since $T_1$ commutes with $S_2,E_2$.
For the two nonadjacent fork vertices, the images of $r_1,r_2$ commute,
and
\[
 (T_1E_1T_1)E_1=0=E_1(T_1E_1T_1).
\]
Furthermore,
\[
 (T_1T_2S_1)E_1=\epsilon E_1
   =E_1(T_1T_2S_1),\qquad
 S_1(T_1E_1T_1)=\epsilon T_1E_1T_1
   =(T_1E_1T_1)S_1.
\]
All remaining nonadjacent relations involve disjoint positions. Thus $q_\epsilon$ preserves all defining relations of $A$, including $r_{\kappa}e_{\kappa}=e_{\kappa}=e_{\kappa}r_{\kappa}$ and
$e_{\kappa}^2=\delta e_{\kappa}$. By Lemma~\ref{lem:D-colored-generators}, we see that $q_\epsilon$ is surjective.
\end{proof}

\subsubsection{The two quotients and the ordinary Brauer ideal}

We compare the colored models with the decorated diagram basis of type
$\mathscr{D}$, using the following decoration convention.

\begin{definition}\label{def:decorated-matchings}
\cite[Definition~4.2]{CGWDiagrams2007}
A
\textbf{decorated Brauer matching} assigns a label $0$ or $1$ to
each strand. It is \textbf{evenly decorated} when the number of
strands labelled $1$ is even. Write $T_n$ for the set of such
evenly decorated matchings on $n$ strands.
\end{definition}

Let $\Theta=Ae_1e_2A$ and let $J$ be the ideal generated by all $e_{\kappa}$.
We use two precise consequences of the Cohen--Gijsbers--Wales diagram model
\cite[Definitions~4.2, 4.4 and Proposition~4.8]{CGWDiagrams2007}.
The algebra $A$ has a vector-space basis consisting of
\[
 T_n\ \sqcup\ \xi T_n^{=}\ \sqcup\ \theta U_n^{=},
\]
where $T_n$ is as in Definition~\ref{def:decorated-matchings},
$T_n^{=}$ consists of those with a horizontal strand, and $U_n^{=}$
consists of the undecorated horizontal matchings. Their cardinalities are
\[
 |T_n|=2^{n-1}b_n,\quad |T_n^{=}|=2^{n-1}(b_n-n!),
 \quad |U_n^{=}|=b_n-n!.
\]
Here $\theta U_n^{=}$ spans $\Theta$. Multiplication by the formal central
symbol $\xi$ is a well-defined bimodule endomorphism of $J$, and
\[
 \xi^2=\delta^2\quad\hbox{on }J/\Theta,
 \qquad \xi e_2=\delta r_1e_2.
\]
The latter identity is the mate-reflection identity in the diagram model;
it also follows from the explicit formula immediately before Lemma~4.7
of the cited paper. The symbol $\xi$ need not itself be an element of $A$.
Only its action on $J$ is used here.

Since the Dynkin diagram is connected, $J=Ae_2A$. Lemma~\ref{lem:D-colored-map}
therefore gives, for every $x\in J$,
\[
 q_\epsilon(\xi x)=\epsilon\delta q_\epsilon(x).
\]
In $J/\Theta$, on each pair of basis vectors $\{t,\xi t\}$, the operator $\xi$ has eigenvalues $\delta$ and $-\delta$, so its two eigenspaces have the same dimension
$h=2^{n-1}(b_n-n!)$. Since $q_\epsilon$ is surjective,
the kernel of the induced map on $A/\Theta$ has dimension
\[
 \dim(A/\Theta)-\dim F_\epsilon=h,
\]
and hence coincides with the $-\epsilon\delta$-eigenspace of $\xi$ on $J/\Theta$.
Thus the two kernels intersect trivially in $A/\Theta$.

Modulo their horizontal ideals, both $F_+$ and $F_-$ have the same quotient $\mathbb C[W(\mathscr{D}_n)]$, and $q_+$ and $q_-$ induce the same map onto this quotient. Therefore, $(q_+,q_-)$ induces an injective homomorphism from $A/\Theta$ into the fiber product. Since both $A/\Theta$ and the fiber product have dimension
$2^n b_n-2^{n-1}n!$, we get the isomorphism
\begin{equation}\label{eq:D-fiber-product}
 A/\Theta\ \cong\quad
 F_+\mathop{\times}_{\mathbb C[W(\mathscr{D}_n)]}F_-.
\end{equation}

By \cite[Proposition 6.4(ii)]{CGW2007}, if $I_n(\delta)$ is the horizontal
ideal of the ordinary Brauer algebra $B_n(\delta)$, then
\begin{equation}\label{eq:D-theta-ordinary}
 \Theta\cong I_n(\delta),\qquad
 u(\delta^{-1}e_1)e_2v\longmapsto u e_2v.
\end{equation}
In particular, the parameter on the right is $\delta$, not $\delta^2$.

\subsubsection{Semisimplicity of the factors}

Let us recall three basic facts about semisimplicity. First, as shown in the proof of Lemma~\ref{lem:sandwich-layers}, a finite-dimensional unital algebra $A$ is semisimple if and only if both a two-sided ideal $I$ and the quotient $A/I$ are semisimple. Second, a fiber
product of two algebras over a common semisimple quotient is semisimple if and only if
both factors are semisimple. Third, for a finite group $G=\{1,\tau\}\cong C_2$ acting on a finite-dimensional complex algebra $B$, the crossed product $B\rtimes C_2$
is semisimple if and only if $B$ is semisimple: if $B$ is semisimple, we average over $C_2$ as in Maschke's theorem, whereas if $\operatorname{rad}(B)\ne0$, then $\operatorname{rad}(B)\rtimes C_2$
is a nonzero nilpotent ideal of $B\rtimes C_2$.

\begin{proposition}\label{prop:D-three-factors}
Let $n\ge4$, $\delta\in\mathbb C^\times$, and
$A=\Br(\mathscr{D}_n;\delta)$. Then
\begin{equation}\label{eq:D-three-factors}
 A\text{ is semisimple}\quad\Longleftrightarrow\quad
 B_n(\delta),\ C_+(n,\delta/2),\text{ and }C_-(n,\delta/2)
 \text{ are semisimple}.
\end{equation}
\end{proposition}
\begin{proof}
By \eqref{eq:D-theta-ordinary}, the ideal $\Theta$ is semisimple if and
only if $B_n(\delta)$ is semisimple, since the quotient of $B_n(\delta)$
by its horizontal ideal is the semisimple group algebra
$\mathbb C[\mathfrak S_n]$. By \eqref{eq:D-fiber-product}, $A/\Theta$
is semisimple if and only if both $F_+$ and $F_-$ are semisimple.
By Lemma~\ref{lem:D-colored-generators} and the crossed-product property above, $F_+$ and $F_-$ are semisimple if and only if $C_+(n,\delta/2)$ and
$C_-(n,\delta/2)$ are semisimple, respectively. Combining these equivalences, we get our result.
\end{proof}

We next apply the ordinary criterion and the walled Brauer criterion of
Cox, De Visscher, Doty, and Martin \cite[page~197, Theorem~6.3]{CDDM2008} to the two
colored factors.
\begin{proposition}\label{prop:colored-semisimplicity}
Let $d\in\mathbb C^\times$. For $n\ge2$,
\[
 C_+(n,d)\text{ is semisimple}\quad\Longleftrightarrow\quad d\notin Z(n).
\]
For $n\ge4$,
\[
 C_-(n,d)\text{ is semisimple}\quad\Longleftrightarrow\quad
 d\notin\{2-n,3-n,\ldots,n-2\}.
\]
\end{proposition}
\begin{proof}
In $C_+(n,d)$, each horizontal arc joins equal colors. By the edge Fourier
change in Lemma~\ref{lem:fourier-gram}, its cell pairings are the matrices
$H_{n,j,\theta}(2d,0)$, with both ordinary parameters equal to $d$.
If $d\notin Z(n)$, then $d\notin Z(v)$ for $2\le v\le n$, so every
ordinary factor is invertible by Theorem~\ref{thm:typeA-rui-si}.
Conversely, if $d\in Z(n)$, we choose a singular ordinary matrix
$D_{n,a,\lambda}(d)$. Taking $j=a$ and $\theta=(\lambda,\varnothing)$, we get
the singular block $2^aD_{n,a,\lambda}(d)$.

In $C_-(n,d)$, each horizontal strand joins opposite colors, so the number of
$+1$ colors in the top word is equal to that in the bottom word. Since all words
with $r$ colors $+1$ are isomorphic via permutation diagrams and the
endomorphism algebra of $+^r-^{n-r}$ is $B_{r,n-r}^{\rm wall}(d)$, we obtain
\begin{equation}\label{eq:D-walled-decomposition}
 C_-(n,d)\cong\bigoplus_{r=0}^n
 \operatorname{Mat}_{\binom nr}\bigl(B_{r,n-r}^{\rm wall}(d)\bigr).
\end{equation}
The summands with $r=0,n$ are group algebras. For $1\le r\le n-1$ and $d\ne0$,
by \cite[page~197, Theorem~6.3]{CDDM2008}, $B_{r,n-r}^{\rm wall}(d)$ is semisimple if and only if
$d\notin\mathbb Z$ or $|d|>n-2$. Therefore, $C_-(n,d)$ is semisimple
if and only if $d\notin\{2-n,3-n,\ldots,n-2\}$.
\end{proof}

For the type $\mathscr{D}$ criterion, the forbidden parameters from $C_-$ are
already accounted for by the other two factors. Indeed, if an integer $k$ satisfies $2-n\leq k\leq n-2$ and
$k\notin Z(n)$, then $k$ is odd and $k\leq3-n$. Therefore
$2k$ is even and $4-2n\leq2k\leq n-2$, so $2k\in Z(n)$. Hence
\[
 2\{2-n,\ldots,n-2\}\subseteq Z(n)\cup2Z(n).
\]

\begin{theorem}\label{thm:D-complete}
For $n\geq4$ and $\delta\in\mathbb C^\times$, we have
\[
 \operatorname{Br}(\mathscr{D}_n;\delta)\text{ is semisimple}
 \quad\Longleftrightarrow\quad
 \delta\notin Z(n)\cup2Z(n).
\]
\end{theorem}
\begin{proof}
If $A=\operatorname{Br}(\mathscr{D}_n;\delta)$ is semisimple, then by \eqref{eq:D-three-factors}, both $B_n(\delta)$ and $C_+(n,\delta/2)$ are semisimple, which gives $\delta\notin Z(n)$ and $\delta/2\notin Z(n)$, namely $\delta\notin Z(n)\cup2Z(n)$.
Conversely, if $\delta\notin Z(n)\cup2Z(n)$, then $B_n(\delta)$ and $C_+(n,\delta/2)$ are semisimple. Moreover, since $2\{2-n,\ldots,n-2\}\subseteq Z(n)\cup2Z(n)$, we have $\delta/2\notin\{2-n,\ldots,n-2\}$, so $C_-(n,\delta/2)$ is also semisimple by Proposition~\ref{prop:colored-semisimplicity}. Applying
\eqref{eq:D-three-factors}, we get our result.
\end{proof}

\begin{remark}
The ordinary ideal in \eqref{eq:D-theta-ordinary} has parameter
$\delta$. For example, the maximal-frame form in type $\mathscr{D}_4$ is
\[
 \delta^3\begin{pmatrix}\delta&1&1\\1&\delta&1\\1&1&\delta\end{pmatrix}
 =\delta^2\begin{pmatrix}\delta^2&\delta&\delta\\
 \delta&\delta^2&\delta\\\delta&\delta&\delta^2\end{pmatrix}.
\]
After removing $\delta^2$, this is the ordinary $B_4(\delta)$
zero-through form. The published statement
\cite[Theorem~5.1(i)]{LiuMorita2025} uses the condition $\delta^2\in Z(n)$.
We do not invoke that assertion here: the ideal isomorphism
\eqref{eq:D-theta-ordinary} identifies the ordinary Brauer parameter as
$\delta$.
\end{remark}

\section{Exceptional and dihedral types}
\label{sec:exceptional-cases}
\label{subsec:type-E}
\label{app:exceptional}

In this section, we determine semisimplicity for the exceptional types using admissible
root-set layers. We first use central operators and perpendicular contractions to restrict
the candidate parameters, and then settle the remaining cases by exact matrix calculations. For the dihedral models, we compute the
layer pairings by counting loops.

Throughout the simply laced and folded exceptional cases in this section, we assume
$\delta\in\mathbb C^\times$ and use the presentations and normal-form bases from the cited papers. For the dihedral diagram models, we will treat all complex parameters $\delta\in\mathbb C$.

\subsection{Simply laced exceptional types: layer modules}
Assume $\delta\in\mathbb C^\times$. We normalize the roots of a simply laced root system $\Phi$ by $(a,a)=2$, and we identify $a$ with $-a$
when referring to root lines. Write $W=W(\Phi)$ and $A=\Br(\Phi;\delta)$.
For $a=w\alpha_i$, we set
$r_a=wr_iw^{-1}$ and $e_a=we_iw^{-1}$; by \cite[equation~(4) and Lemma~4.2]{CFW2009}, these elements are independent
of the choice of $w$ and $i$.
Moreover, by \cite[Lemma~4.3(iii)]{CFW2009}, $e_a$ and $e_b$ commute whenever $a$ and $b$ are orthogonal.
\begin{definition}\label{def:admissible-root-sets}
\cite[Section~1 and Lemma~2.1]{CFW2009}
An orthogonal
set $X$ of positive roots is \textbf{admissible} if, whenever
$b_1,b_2,b_3\in X$ are distinct and $a\in\Phi$ satisfies
$|(a,b_j)|=1$ for all three indices, $X$ also contains the positive
representative of
$r_a r_{b_1}r_{b_2}r_{b_3}a$.
The \textbf{admissible closure} $X^{\mathrm{cl}}$ is the smallest
admissible set containing $X$ \cite[Definition~2.2]{CFW2009}.
We also use \textbf{$k$-frame} for a set of $k$ pairwise orthogonal
root lines; a full frame spans the ambient root space.
\end{definition}

The admissible sets index the boundary data in the CFW normal forms.
Their root products give the idempotents used to detect specialized layers.
Set
\[
 e_X=\prod_{a\in X}e_a,\qquad
 f_X=\delta^{-|X|}e_X\quad(\delta\ne0).
\]
Thus $e_X^2=\delta^{|X|}e_X$ and $f_X^2=f_X$.
\begin{definition}\label{def:coefficient-group}
For admissible $X$, define its \textbf{Weyl annihilator}
$A_X$ in the unspecialized algebra over $\mathbb C[z,z^{-1}]$,
before setting $z=\delta$. Its \textbf{coefficient group} is the quotient
$H_X$ in
\[
 A_X=\{w\in W:we_X=e_X\},\qquad
 H_X=N_W(X)/A_X.
\]
Here $N_W(X)=\{w\in W:wX=X\}$ is the setwise stabilizer of the
root lines. This is \cite[equation~(7) and Proposition~4.7]{CFW2009}.
\end{definition}

By Proposition~4.7 of \cite{CFW2009}, $A_X$ is a normal subgroup of
$N_W(X)$ and the quotient $H_X=N_W(X)/A_X$ is isomorphic to the coefficient Weyl group
(note that $A_X$ need not coincide with
the pointwise stabilizer of $X$).

We use the normal-form basis in
\cite[Proposition~4.9 and Corollary~5.5]{CFW2009} and the layer data in
\cite[Lemma~1.3 and Table~3]{CFW2009} for $\mathscr{E}_6,\mathscr{E}_7,\mathscr{E}_8$.
The proof of inequivalence in \cite[Proposition~5.4]{CFW2009}
was corrected in \cite{CFWerr2019}; the basis and layer data are unchanged.
For each orbit $\mathcal B=WX$, we choose a representative $w_Y\in W$ with $w_YX=Y$ for each $Y\in\mathcal B$.
Each normal form has two boundary labels in $\mathcal B$ and a
middle element in $H_X$, and the span of the orbit $\mathcal B$ modulo the strata of strictly larger
cardinality is an inflation of $\mathbb C[H_X]$.
If we choose different representatives $w_Y$, the fiber bases change by invertible transformations, which preserves the rank of each layer pairing. We list the complete nonempty layer data in the following table:
\[
\begin{array}{c|r|r|c}
\text{type}&|X|&|\mathcal B|&H_X\\\hline
 \mathscr{E}_6&1&36&\mathfrak S_6\\
    &2&270&\mathfrak S_3\\
    &4&135&1\\\hline
 \mathscr{E}_7&1&63&W(\mathscr{D}_6)\\
    &2&945&\mathfrak S_4\times C_2\\
    &3&315&\mathfrak S_3\\
    &4&945&C_2\\
    &7&135&1\\\hline
 \mathscr{E}_8&1&120&W(\mathscr{E}_7)\\
    &2&3780&\mathfrak S_6\\
    &4&9450&\mathfrak S_3\\
    &8&2025&1.
\end{array}
\]
Together with the empty-set layer $\mathbb C[W]$, we have the dimension formula
\[
 \dim A=|W|+\sum_{\mathcal B\ne\{\varnothing\}}
       |\mathcal B|^2|H_X|,
\]
which confirms that these orbits exhaust the entire algebra $A$.

The CFW normal-form action gives the layer representations. For an irreducible representation $\rho_\theta$
of $H_X$, the Weyl group $W$ acts by
\[
 g(v_Y\otimes u)=v_{gY}\otimes
 \rho_\theta\!\left(\overline{w_{gY}^{-1}gw_Y}\right)u,
\]
where the bar denotes the quotient map $N_W(X)\to H_X$.
On the fiber at $Y$, a root contraction $e_a$
acts as $\delta$ if $a\in Y$, as zero in this layer if $a\perp Y$,
and otherwise by the Weyl transport $r_b r_a$, where we choose
$b\in Y$ with $(a,b)\ne0$. By the CFW normal-form relations, this
action is independent of the choice of $b$. Letting $\theta$ run over $\Irr(H_X)$,
we obtain all coefficient blocks for the layer $\mathcal B$.

Fix a nonempty orbit $\mathcal B=WX$ and put $t=|X|$.
For $\theta\in\Irr(H_X)$, the corresponding layer module is
\[
 V_{\mathcal B,\theta}=\bigoplus_{Y\in\mathcal B}V_Y,
 \qquad V_Y\simeq U_\theta.
\]
The product $e_Y=\prod_{a\in Y}e_a$ has image $V_Y$, and restricts to
$\delta^tI$ there. Set
\[
 T_{\mathcal B,\theta}(\delta)=\sum_{Y\in\mathcal B}e_Y,
 \qquad K_{\mathcal B,\theta}=\bigcap_Y\ker e_Y.
\]
\begin{proposition}\label{prop:exceptional-layer-test}
Let $\Phi$ have type $\mathscr{E}_6$, $\mathscr{E}_7$, or $\mathscr{E}_8$, and let
$\delta\in\mathbb C^\times$. With the layer modules and operators above,
\begin{equation}\label{eq:app-layer-criterion}
 \Br(\Phi;\delta)\text{ is semisimple}
 \quad\Longleftrightarrow\quad
 \det T_{\mathcal B,\theta}(\delta)\ne0
 \text{ for all }(\mathcal B,\theta).
\end{equation}
\end{proposition}
\begin{proof}
Since the images of the summands $e_Y$ lie in distinct fibers $V_Y$, we have
$K_{\mathcal B,\theta}=\ker T_{\mathcal B,\theta}$.
By Lemma~\ref{lem:specialized-layer-detection} below,
$K_{\mathcal B,\theta}$ is an $A$-submodule of $V_{\mathcal B,\theta}$. If $K_{\mathcal B,\theta}\ne0$, then it has no
$A$-invariant complement in $V_{\mathcal B,\theta}$, because any invariant complement must contain
$e_YV_{\mathcal B,\theta}=V_Y$ for every $Y\in\mathcal B$ and hence equal $V_{\mathcal B,\theta}$.
Conversely, suppose $K_{\mathcal B,\theta}=0$, and let $U$ be a nonzero submodule of $V_{\mathcal B,\theta}$.
We choose a nonzero vector $u\in U$. Since $\bigcap_Y\ker e_Y=0$, we have $e_Yu\ne0$ for some $Y\in\mathcal B$, and $e_Yu\in U\cap V_Y$. The coefficient-group action transported to $V_Y$ is irreducible, so
we get $V_Y\subseteq U$, and by Weyl transport we see that $U$ contains every fiber. Thus $V_{\mathcal B,\theta}$ is simple.
For a fixed orbit, applying the idempotent $f_X=\delta^{-t}e_X$ to a simple layer module recovers
the coefficient module $U_\theta$, so distinct
$\theta\in\Irr(H_X)$ give nonisomorphic simple modules, while modules from different
strata are distinguished by the ideal filtration in
Lemma~\ref{lem:specialized-layer-detection}. Therefore, if every $T_{\mathcal B,\theta}(\delta)$ is
invertible, the squared dimensions of these simple modules together with those of the group
quotient sum to $\dim\Br(\Phi;\delta)$, so $\Br(\Phi;\delta)$ is semisimple.
If some $T_{\mathcal B,\theta}(\delta)$ is singular, then $K_{\mathcal B,\theta}$ is a nonzero submodule with no invariant complement, so $\Br(\Phi;\delta)$ is not semisimple.
\end{proof}

\begin{lemma}
\label{lem:specialized-layer-detection}
Fix $\delta\ne0$ in one of the simply laced exceptional types.
Let $J_{\ge t}$ be the span of the normal forms whose admissible
boundary sets have cardinality at least $t$, and put
$J_{>t}=J_{\ge t+1}$. These are two-sided ideals. For an admissible
set $X$ of size $t$, write $f_X=\delta^{-t}e_X$.
In the quotient by $J_{>t}$, the ideal $J_{\ge t}/J_{>t}$ is generated
by these idempotents; for a representative $X$ of a size-$t$ orbit,
its coefficient corner is $f_X(A/J_{>t})f_X\cong\mathbb C[H_X]$.

Every simple $A$-module has a unique stratum cardinality $s$:
$J_{>s}L=0$ and $f_XL\ne0$ for a set $X$ of size $s$.
The nonzero corner module $f_XL$ is simple over $\mathbb C[H_X]$.
For a layer module $V=V_{\mathcal B,\theta}$ of cardinality $t$, its
submodule $K=\bigcap_{Y\in\mathcal B}\ker e_Y$ is annihilated by
$J_{\ge t}$. Consequently every composition factor of $K$ has
cardinality strictly smaller than $t$.
\end{lemma}
\begin{proof}
Since multiplying two normal forms never decreases the cardinality of the admissible boundary set, $J_{\ge t}$ and $J_{>t}$ are two-sided ideals. Modulo $J_{>t}$,
the normal forms in a fixed orbit have the form $u e_X h v$, so
their multiplication gives $f_X(A/J_{>t})f_X\cong\mathbb C[H_X]$ and shows
that the conjugate idempotents $f_Y$ generate $J_{\ge t}/J_{>t}$.
Since we only divide by powers of $\delta\ne0$, these properties hold at every nonzero complex specialization.

For a simple $A$-module $L$, let $s$ be the largest cardinality
such that $J_{\ge s}L\ne0$ (with $J_{\ge0}=A$).
Then $J_{>s}L=0$, and since $J_{\ge s}/J_{>s}$ is generated by the idempotents $f_X$ with $|X|=s$,
we have $f_XL\ne0$ for some such $X$. Moreover, if
$0\ne U\subseteq f_XL$ is a submodule over $f_X(A/J_{>s})f_X$, then the simplicity of $L$ gives
$(A/J_{>s})U=L$, and hence
$f_XL=f_X(A/J_{>s})f_XU=U$. Thus $f_XL$ is a simple $\mathbb C[H_X]$-module.

By definition of the layer action, $J_{>t}$ annihilates $V$. Any normal form
of cardinality $t$ can be written as $u e_X h v$ with $u,h,v\in W$. Conjugating the root product, we have
\[
 u e_X h v=u h v\,e_{(hv)^{-1}X}.
\]
Since $(hv)^{-1}X\in\mathcal B$, this element
annihilates $K$. In types $\mathscr{E}_6,\mathscr{E}_7,$ and $\mathscr{E}_8$, there is only one admissible-set
orbit at each allowed positive cardinality, so we get $J_{\ge t}K=0$.
Since each $e_Y$ belongs to $J_{\ge t}$, we see that
$K=\operatorname{Ann}_{J_{\ge t}}V$ is an $A$-submodule of $V$.
By the characterization of $s$ above, every composition factor of $K$ has
stratum cardinality strictly smaller than $t$.
\end{proof}

\begin{remark}
For any idempotent $f$, the functor $M\mapsto fM$ is exact: a lift of
$fm$ through a surjection is obtained by applying $f$ to any lift of $m$.
It follows that $fM=0$ if and only if $f$ annihilates every composition
factor of $M$. We use this exactness property in our perpendicular-contraction reduction: in particular, if every composition factor of $K$
has stratum cardinality zero, then every idempotent $\delta^{-1}e_a$ annihilates $K$ itself.
\end{remark}

\begin{lemma}
\label{lem:central-collision}
For a simply laced finite root system, set
$\Omega=\sum_{a\in\Phi^+}(r_a-e_a)$, indexing positive roots as root lines.
Then $\Omega$ is central. On the layer $V_{\mathcal B,\theta}$ of
admissible cardinality $t$, it acts as
$t(1-\delta)+c_{t,\theta}$, where $c_{t,\theta}$ is the scalar of the
perpendicular reflection sum on the coefficient representation.
\end{lemma}
\begin{proof}
Since conjugation by $W$ permutes the positive root lines, $\Omega$ commutes with every
reflection. Now fix a root $a$. The summands indexed by $a$ and by roots
orthogonal to $a$ commute with $e_a$. The remaining root lines split
into pairs $\{b,r_a b\}$. In the rank-two subsystem of type $\mathscr{A}_2$ spanned by $a$ and $b$, the
Brauer relations give
\[
 e_a(r_b+r_{r_a b}-e_b-e_{r_a b})=0
 =(r_b+r_{r_a b}-e_b-e_{r_a b})e_a.
\]
Indeed, for adjacent simple generators $r_1,r_2,e_1,e_2$, we have
$e_1r_2r_1=e_1e_2$ and $e_1e_2r_1=e_1r_2$, while
$r_{r_1b}=r_1r_2r_1$ and $e_{r_1b}=r_1e_2r_1$.
These identities give the cancellation with $e_a$ on the left; applying
the anti-involution $*$ gives the cancellation with $e_a$ on the right. Thus $\Omega$ also commutes with every
$e_a$, and hence is central.

For an admissible set $X=\{a_1,\ldots,a_t\}$, we partition the root lines
outside $X\cup X^\perp$ according to the first $a_j$ to which they
are not orthogonal. Since $a_1,\ldots,a_t$ are mutually orthogonal, each part is invariant under $r_{a_j}$, so the same rank-two cancellation shows
that all these paired summands annihilate $e_X$.
Meanwhile, the roots in $X$ contribute $t(1-\delta)e_X$.
For any root $b\in X^\perp\setminus X$, the product
$e_b e_X$ lies in a strictly larger admissible stratum and hence
acts as zero on $V_{\mathcal B,\theta}$, whereas the sum of the reflections $r_b$ over positive root lines
$b\in X^\perp\setminus X$ acts on the coefficient fiber $V_X$ by the scalar $c_{t,\theta}$. Therefore, $\Omega$ acts as $t(1-\delta)+c_{t,\theta}$ on $V_X$, and by centrality and Weyl transport it acts by the same scalar on all of $V_{\mathcal B,\theta}$.
\end{proof}

Combining Lemmas~\ref{lem:specialized-layer-detection} and~\ref{lem:central-collision}, we see that if $K_{\mathcal B,\theta}$ has a composition factor of stratum cardinality $s<t$ with coefficient label $\eta$, then $\delta$ must satisfy the central collision condition
\begin{equation}\label{eq:app-collision}
 \delta=1+\frac{c_{t,\theta}-c_{s,\eta}}{t-s}.
\end{equation}
Since the layer operator on the regular representation of $H_X$ is a monic polynomial in $\delta$ with integer matrix coefficients, every nonsemisimple parameter $\delta$ is an
algebraic integer; since \eqref{eq:app-collision} shows that $\delta\in\mathbb Q$, any such parameter must be a rational
integer.

\begin{proposition}\label{prop:app-contraction}
Let $K=K_{\mathcal B,\theta}$ be the kernel of a layer module $V$, let $X$ be an orthogonal root set of size $k$,
and put $P=\Phi\cap X^\perp$. If $e_XK\ne0$, then the perpendicular Brauer algebra $\Br(P;\delta)$ is
not semisimple.
\end{proposition}
\begin{proof}
Since the generators of $\Br(P;\delta)$ commute with $e_X$, the subspace
\[
 M=e_XV=\bigoplus_{B\supseteq X}V_B
\]
and its subspace $K'=e_XK$ are $\Br(P;\delta)$-modules.
For each $B\in\mathcal B$ with $B\supseteq X$, the product $e_{B\setminus X}$ lies in the
image of $\Br(P;\delta)$ and satisfies
\[
 e_{B\setminus X}|_M=\delta^{-k}e_B|_M,
 \qquad e_{B\setminus X}M=V_B,
 \qquad e_{B\setminus X}K'=0.
\]
If $K'$ had a $\Br(P;\delta)$-invariant complement $U$ in $M=K'\oplus U$,
then every $V_B=e_{B\setminus X}M=e_{B\setminus X}U$ would lie in $U$, which implies $M=U$ and $K'=0$, a contradiction.
Thus $M$ is not semisimple as a $\Br(P;\delta)$-module, so $\Br(P;\delta)$ is not semisimple.
\end{proof}

For type $\mathscr{D}$, we have
\begin{equation}\label{eq:app-D-input}
 \Br(\mathscr{D}_n;\delta)\text{ is semisimple}
 \quad\Longleftrightarrow\quad
 \delta\notin Z(n)\cup2Z(n),\qquad n\ge4,
\end{equation}
where
\[
 Z(n)=\{4-2n,\ldots,n-2\}
 \setminus\{z\text{ odd}:4-2n<z\le3-n\}.
\]
This is Theorem~\ref{thm:D-complete}.
In particular,
\[
 \operatorname{Bad}(\mathscr{D}_6)=
 \{-16,-12,-8,-6,-4,-2,-1,1,2,3,4,6,8\}.
\]

\subsubsection{Exact certificates for finite calculations}
\label{subsec:exact-certificate-principles}
We compute all finite ranks below over specified exact fields. The following
lemma gives the rank bounds and polynomial identities used to draw
conclusions in characteristic zero.

\begin{lemma}
\label{lem:exact-matrix-certificates}
Let $M\in\operatorname{Mat}_N(\mathbb Z)$ and let $p$ be prime.
\begin{enumerate}
\item If $M$ has rank $N$ modulo $p$, then $M$ is invertible over
$\mathbb Q$ and $\mathbb C$.
\item If $M$ has rank at least $r$ modulo $p$, and $N-r$ linearly
independent rational vectors are verified to belong to $\ker M$,
then $\operatorname{rank}_{\mathbb Q}M=r$.
\item Suppose $M(z)\in\operatorname{Mat}_N(\mathbb Q[z])$ and
$F(z)\in\mathbb Q[z]$ satisfy $\deg\det M(z)\le D$ and
$\deg F(z)\le D$. If
$\det M(a)=F(a)$ is verified exactly at $D+1$ distinct rational numbers,
then $\det M(z)=F(z)$. In particular, entry degree at most
$d$ gives the bound $D=Nd$.
\end{enumerate}
\end{lemma}
\begin{proof}
Since any minor of $M$ that is nonzero modulo $p$ is also nonzero over $\mathbb Z$, we get the
rank lower bound and (1). In (2), the $N-r$ linearly independent kernel vectors
give $\operatorname{rank}_{\mathbb Q}M\le r$, which yields $\operatorname{rank}_{\mathbb Q}M=r$.
In (3), the polynomial $\det M(z)-F(z)$ has degree at most $D$ and vanishes at
$D+1$ distinct points, so it is identically zero.
\end{proof}

For a polynomial matrix identity whose entries have degree at most
$d$, exact verification at $d+1$ distinct rational points suffices.
This uses the entry-degree bound, rather than the determinant-degree bound.

To bound the root multiplicities of $\det M(z)$, we use truncated local jet
matrices. Suppose $\det M(z)$ is not identically zero and let $a\in\mathbb Q$.
If the local Smith invariant factors of $M(z)$ at $z=a$ are
$(z-a)^{b_1},\ldots,(z-a)^{b_N}$, then the nullity of multiplication
by $M(z)$ on $(\mathbb Q[z]/(z-a)^k)^N$ is equal to
\[
 \sum_{j=1}^N\min\{k,b_j\},
\]
because the local Smith transformations remain invertible modulo $(z-a)^k$. In particular, this nullity is a lower bound for
$\operatorname{ord}_{z=a}\det M(z)=\sum_{j=1}^N b_j$. When $M(z)$ has monic determinant and the sum of
these lower bounds over a set of roots equals $\deg\det M(z)$, the resulting factorization of $\det M(z)$ is exact.

In our branching calculations below, we use the explicit geometric subgroup embeddings $H\le W$. Given a character $\theta$ of $H$, the multiplicity
of $\chi\in\operatorname{Irr}(W)$ in $\operatorname{Ind}_H^W\theta$ is given by
\[
 \frac1{|H|}\sum_{h\in H}\chi(h)\theta(h^{-1}).
\]
To evaluate this sum over conjugacy classes of $H$, we determine the conjugacy class of each $h\in H$ in $W$ directly from the root permutation action. Specifically, in type $\mathscr{E}_8$, we act with the class representatives on the $240$ roots; since all $112$ conjugacy classes of $W(\mathscr{E}_8)$ have distinct cycle types on the $240$ roots, the cycle type uniquely identifies the target conjugacy class in $W(\mathscr{E}_8)$.

We apply Lemma~\ref{lem:exact-matrix-certificates}(1) to the three
residual $\mathscr{E}_7$ matrices at $z=-7$. For $\mathscr{F}_4$, we use the exact integer
evaluations in part~(3). The $\mathscr{E}_6$ and $\mathscr{H}_4$ calculations instead use
the rank, jet, and spectral certificates described below.

\subsubsection{Singleton spectra for \texorpdfstring{$\mathscr{E}_6,\mathscr{E}_7,\mathscr{E}_8$}{E6, E7, E8}}
\label{subsec:singleton-spectra}
Fix a root line $a$, let $s_0=r_a$, and put $H=W(\Phi\cap a^\perp)$.
The singleton module satisfies
\[
 V_\theta|_W=\Ind_{\langle s_0\rangle\times H}^{W}
       (1\otimes U_\theta).
\]
Write $b_{\lambda\theta}$ for the induced multiplicity, and $c_\lambda$
and $c_\theta$ for the reflection-sum scalars in $W$ and $H$.
By the singleton fiber action, we have
\[
 T_\theta=\sum_{a\in\Phi^+}r_a+(\delta-1-c_\theta)I,
\]
and therefore
\begin{equation}\label{eq:app-singleton-det}
 \det T_\theta(\delta)=
 \prod_{\lambda\in\Irr(W)}
   (\delta-1-c_\theta+c_\lambda)^{b_{\lambda\theta}\dim U_\lambda}.
\end{equation}
We compute the induction multiplicities $b_{\lambda\theta}$ by the formula
\[
 b_{\lambda\theta}=
 \frac1{2|H|}\sum_{C\subset H}|C|
  \bigl(\chi_\lambda(C)+\chi_\lambda(s_0C)\bigr)\chi_\theta(C),
\]
where we determine class fusion from the explicit root permutations
and take the character tables from the $q=1$ specialization of CHEVIE
\cite{CHEVIE}. We check character orthogonality, nonnegative integral
multiplicities, and the dimension of each induced representation.
In type $\mathscr{E}_7$, whenever two irreducible characters of $W(\mathscr{D}_6)$ are paired together, they share the same reflection scalar $c_\theta$, so grouping them preserves the set of roots of \eqref{eq:app-singleton-det}. In type $\mathscr{E}_8$, all $112$ conjugacy classes are distinguished by their cycle types on the $240$ roots, from which we obtain the $112\times60$ singleton branching table.

Taking the union of all nonzero roots of \eqref{eq:app-singleton-det} over $\theta\in\Irr(H)$, we obtain the following sets of parameters:
\begin{align*}
 S_6=\{&-20,-14,-12,-8,-6,-5,-4,-3,-2,1,2,4,5,6,7,8,10\},\\
 S_7=\{&-32,-24,-20,-16,-14,-12,-10,-8,-6,-5,-4,-3,-2,-1,\\
       &1,2,3,4,5,6,7,8,10,12,14,16\},\\
 S_8=\{&-56,-44,-36,-32,-26,-24,-20,-18,-16,-14,-12,-11,\\
       &-10,-9,-8,-7,-6,-5,-4,-2,-1,1,2,3,4,6,7,8,9,10,\\
       &11,12,13,14,16,18,20,22,26,28\},
\end{align*}
which have cardinalities $17,26,$ and $40$, respectively. By \eqref{eq:app-layer-criterion}, $\Br(\mathscr{E}_\kappa;\delta)$ is not semisimple for any $\delta\in S_\kappa$.

\begin{theorem}\label{thm:app-exceptional}
For $\kappa\in\{6,7,8\}$ and $\delta\in\C^*$, we have
\[
 \Br(\mathscr{E}_{\kappa};\delta)\text{ is semisimple}
 \quad\Longleftrightarrow\quad \delta\notin S_{\kappa}.
\]
\end{theorem}
In the next three subsections, we prove that all higher-cardinality layers of $\mathscr{E}_6,\mathscr{E}_7,$ and $\mathscr{E}_8$ are nonsingular whenever $\delta\notin S_\kappa$.

\subsection{Type \texorpdfstring{$\mathscr{E}_6$}{E6}: Fifteen layer determinants}
For type $\mathscr{E}_6$, we calculate the determinants of the four remaining higher-layer matrices directly.
In type $\mathscr{E}_6$, there are three nonempty admissible-set orbits:
\[
\begin{array}{c|r|c|r}
 t&|\mathcal B|&H_X&|\Irr(H_X)|\\\hline
 1&36&S_6&11\\
 2&270&S_3&3\\
 4&135&1&1
\end{array}
\]
Their squared module dimensions, together with the group quotient, give
\[
 51840+36^2\cdot720+270^2\cdot6+135^2=1440585.
\]
For a pair $X$, the perpendicular $\mathscr{A}_3$ subsystem contains three orthogonal root
pairs, and the permutation action on these three pairs realizes the
coefficient quotient $S_3$. Using its trivial, sign, and standard
representations together with Weyl transport and the CFW contraction rule, we obtain integral matrices
$T_{2,\theta}(\delta)=\delta^2I+\delta B_{\theta,1}+B_{\theta,0}$ for the three pair layers, and we check the defining relations and root-conjugation identities
symbolically in $\delta$.

The four higher-layer determinants are
\begin{align*}
 \det T_{2,1}&=(\delta+14)(\delta+8)^{42}(\delta+5)^{20}
 (\delta+2)^{150}(\delta-1)^{24}(\delta-2)^{243}(\delta-4)^{60},\\
 \det T_{2,\mathrm{sgn}}&=(\delta+4)^{81}(\delta+2)^{250}
 (\delta-1)^{24}(\delta-4)^{165}(\delta-7)^{20},\\
 \det T_{2,\mathrm{std}}&=(\delta+8)^{20}(\delta+5)^{15}
 (\delta+4)^{324}(\delta-1)^{466}(\delta-4)^{240}(\delta-7)^{15},\\
 \det T_{4,1}&=\delta^{135}(\delta+8)(\delta+4)^{36}
 (\delta+2)^{120}(\delta-1)^{134}(\delta-2)^{99}(\delta-4)^{15}.
\end{align*}
We verify these four determinant formulas as follows. For $T_{2,1}$,
the nullities at $\delta=-14,-8,-5,-2,1,2,4$ are
$1,42,20,150,24,243,60$, whose sum equals the degree $540$ of $\det T_{2,1}$.
For $T_{2,\mathrm{sgn}}$, the nullities at
$\delta=-4,-2,1,4,7$ are $81,250,24,165,20$.
Their sum is again $540$.
For $T_{2,\mathrm{std}}$, the nullities at $\delta=-8,-5,-4,1,4,7$ are
$20,15,324,406,240,15$. At $\delta=1$, the order-two local jet matrix
\[
 J_2(1)=\begin{pmatrix}T(1)&0\\T'(1)&T(1)\end{pmatrix}
\]
has nullity $466$ (corresponding to $346$ local elementary divisors of exponent one
and $60$ of exponent two). Since the local jet nullity is a lower
bound for the multiplicity of $\delta=1$ in $\det T_{2,\mathrm{std}}$, replacing $406$ by $466$ exhausts
the total degree $1080$.
In the verification package, we compute these ranks and jet nullities
exactly over $\mathbb Q$ with SageMath.
For $T_{4,1}$, the nullities at $\delta=-8,-4,-2,0,1,2,4$
are $1,36,120,135,134,99,15$, whose sum is $540$.
The calculated multiplicity lower bounds exhaust the degree of each
of the four determinant polynomials. Since these polynomials are monic,
we obtain the displayed factorizations.
Every nonzero root of these four determinants lies in $S_6$.
Combining them with the singleton determinants, we prove Theorem~\ref{thm:app-exceptional} for $\mathscr{E}_6$.

\subsection{Type \texorpdfstring{$\mathscr{E}_7$}{E7}: Three residual nonsingularity certificates}
For type $\mathscr{E}_7$, perpendicular contractions and the central collision formula \eqref{eq:app-collision}
reduce the verification outside $S_7$ to three integer matrices at $\delta=-7$.
For any root $a$, its perpendicular root system is of type $\mathscr{D}_6$, and by \eqref{eq:app-D-input} we have
$\operatorname{Bad}(\mathscr{D}_6)\subseteq S_7$.
Hence, for $\delta\notin S_7$, Proposition~\ref{prop:app-contraction} rules out any
composition factor of positive stratum cardinality in a layer kernel $K_{\mathcal B,\theta}$.
Any remaining composition factor must have cardinality $s=0$, so it is a $W(\mathscr{E}_7)$-module that
satisfies both the collision formula \eqref{eq:app-collision} and the stabilizer branching condition.

For a partition $\theta$, we write
$\operatorname{ct}(\theta)=\sum_{(i,j)\in\theta}(j-i)$ for its content sum. For the pair layer ($t=2$), the coefficient group is
$S_4\times C_2$, with
$c_{2,\theta,\varepsilon}=2\operatorname{ct}(\theta)+\varepsilon$.
From the $60\times10$ branching table, we find that the only collision value outside $S_7$ is $\delta=-7$, which occurs for
$([4],-1)$ via the $W(\mathscr{E}_7)$-character $(56,3)$ and for
$([3,1],+1)$ via $(120,4)$.
For the triple and four-root layers ($t=3,4$), the complete nonzero integral
collision values from \eqref{eq:app-collision} with $s=0$ are listed in the following table:
\[
\begin{array}{c|c|r|l}
 t&\theta&c_{t,\theta}&\text{collision values}\\\hline
 3&1&12&-16,-4,4\\
 3&\mathrm{sgn}&-12&-4,8\\
 3&\mathrm{std}&0&-8,-2,8\\
 4&1&3&-14,-7,-5,-2,1,7\\
 4&\mathrm{sgn}&-3&-8,-5,-3,1,7
\end{array}
\]
The triple quotient $S_3$ acts on three perpendicular $\mathscr{D}_4$ four-frames;
the four-root quotient $C_2$ is the product of orientation signs of the
three perpendicular $\mathscr{A}_1$ lines. The stabilizers have orders $9216$ and
$3072$, and the CFW annihilator has order $1536$ in both cases.
All nonzero integral collisions from the seven-root layer lie in $S_7$
even without imposing branching.

At $\delta=-7$, the three remaining candidate layer matrices have the following ranks over $\mathbb F_{1009}$:
\[
\begin{array}{c|c|r|r}
 t&\theta&\text{order}&\rank_{\mathbb F_{1009}}T(-7)\\\hline
 2&([4],-1)&945&945\\
 2&([3,1],+1)&2835&2835\\
 4&1&945&945
\end{array}
\]
Here the first pair model transports the orientation of the perpendicular
$\mathscr{A}_1$ root line, the second transports the augmentation representation on the
four coordinate rays of the perpendicular $\mathscr{D}_4$ system, and the four-root
model acts by monomial matrices on its $945$ fibers. After checking the defining
relations of $\Br(\mathscr{E}_7;\delta)$ on these three models, we compute their ranks modulo $p=1009$ at $\delta=-7$. By Lemma~\ref{lem:exact-matrix-certificates}(1), full rank over $\mathbb F_{1009}$
implies that all three matrices are invertible over $\C$. This rules out $\delta=-7$ and finishes the proof of Theorem~\ref{thm:app-exceptional} for $\mathscr{E}_7$.

\subsection{Type \texorpdfstring{$\mathscr{E}_8$}{E8}: Reduction to six scalar quadratic evaluations}
In type $\mathscr{E}_8$, as in $\mathscr{E}_7$, we use perpendicular contractions to restrict the candidate parameters.
We then evaluate scalar quadratic polynomials on the relevant central eigenspaces to exclude the remaining cases.
The admissible cardinalities are $0,1,2,4,$ and $8$.
For a base orthogonal pair
$X_0=\{\varepsilon_7-\varepsilon_8,\varepsilon_7+\varepsilon_8\}$,
\[
 N_W(X_0)=(W(\mathscr{B}_6)\times W(\mathscr{B}_2))\cap W(\mathscr{D}_8),\qquad
 |N_W(X_0)|=184320,
\]
its orbit has size $3780$, and its coefficient quotient is $S_6$.
The quotient forgets the coordinate signs and the last two coordinates;
its kernel is the CFW subgroup of order $256$.
For $\theta\vdash6$,
$c_{2,\theta}=2\operatorname{ct}(\theta)$.
The exact pair branching table has size $112\times11$.
The four-root orbit has size $9450$, normalizer order $73728$ and
coefficient quotient $S_3$, acting on the three perpendicular $\mathscr{D}_4$
four-frames; its branching table has size $112\times3$.
Its Weyl-factor collision lists are
\[
\begin{array}{c|l}
1&-26,-14,-11,-8,-5,-1,1,4,7\\
\mathrm{std}&-14,-11,-5,-4,-2,1,4,6,7,13\\
\mathrm{sgn}&-8,-5,-2,3,4,13.
\end{array}
\]
All these values belong to $S_8$.

The complete scalar values used in the remaining collision calculations are
\begin{align*}
 C_0&=\{0\}\cup\pm\{6,8,10,12,15,18,20,22,24,30,36,40,42,48,60,72,90,120\},\\
 C_1&=\{0\}\cup\pm\{3,5,7,9,13,15,17,21,27,33,35,45,63\},\\
 C_2&=\{0,\pm6,\pm10,\pm18,\pm30\},\qquad
 C_4=\{0,\pm12\},\qquad C_8=\{0\}.
\end{align*}
Using the stabilizer branching table for $s=0$ and perpendicular
singleton contraction for $(t,s)=(2,1)$, we find that the only collision values outside $S_8$ are as follows:
\[
\begin{array}{c|c|c}
 t&s&\text{remaining candidates}\\\hline
 2&0&-3,5\\
 2&1&\varnothing\\
 4&0&\varnothing\\
 4&1&-3,5\\
 4&2&\varnothing\\
 8&0,1,2,4&\varnothing
\end{array}
\]
In particular, the eight-root layer ($t=8$) gives no candidate parameters outside $S_8$.

At $\delta\in\{-3,5\}$, the perpendicular algebra for an orthogonal pair is
$\Br(\mathscr{D}_6;\delta)$, which is semisimple by \eqref{eq:app-D-input}. By Proposition~\ref{prop:app-contraction}, we can exclude any
composition factor of stratum cardinality at least two.
For the pair layer ($t=2$), contracting by one root gives an $\mathscr{E}_7$ singleton layer
whose $W(\mathscr{D}_6)$-coefficient representation is inflated along $W(\mathscr{D}_6)\to S_6$.
By the $\mathscr{E}_7$ singleton formula \eqref{eq:app-singleton-det}, the eleven corresponding singleton blocks are
nonsingular at $\delta=-3$ and $\delta=5$.
Similarly, for the four-root layer ($t=4$), contracting by one root gives the corresponding
$\mathscr{E}_7$ triple layer, where both coefficient quotients $S_3$ act on the same three four-frames of the common perpendicular $\mathscr{D}_4$ system.
Since $-3$ and $5$ do not appear in the $\mathscr{E}_7$ triple collision table above and $\Br(\mathscr{D}_6;\delta)$ is semisimple at $\delta\in\{-3,5\}$, the $\mathscr{E}_7$ triple layer is nonsingular at both values.
Let $V$ be either original layer module, with kernel $K$, and fix a root $a$.
For every boundary $B$ containing $a$, put $B'=B\setminus\{a\}$. Then
\[
 e_{B'}(e_aK)=e_BK=0,
\]
so $e_aK$ lies in the joint kernel of the corresponding contracted layer.
The singleton and triple layers just considered are nonsingular at both
parameters, hence $e_aK=0$ for every root $a$. By exactness of multiplication
by $\delta^{-1}e_a$ and Lemma~\ref{lem:specialized-layer-detection}, no
cardinality-one composition factor remains. This excludes the
$(t,s)=(2,1)$ and $(4,1)$ cases.

Only six cases remain, each corresponding to a multiplicity-one Weyl
constituent of a pair layer. The corresponding scalar quadratics are
\begin{center}\small
\begin{tabular}{ccllr}
\toprule
$\delta_0$&$\theta$&Weyl label $(\dim,b)$&$q(\delta)$&$q(\delta_0)$\\\midrule
$-3$&$[4,2]$&$(2240,10)$&$(\delta-2)(\delta+4)$&$-5$\\
$-3$&$[3,2,1]$&$(2835,14)$&$(\delta-1)(\delta+8)$&$-20$\\
$-3$&$[3,2,1]$&$(6075,14)$&$(\delta-1)(\delta+8)$&$-20$\\
$5$&$[2,2,1,1]$&$(2240,28)$&$(\delta-8)(\delta+4)$&$-27$\\
$5$&$[3,2,1]$&$(2835,22)$&$(\delta-1)(\delta-8)$&$-12$\\
$5$&$[3,2,1]$&$(6075,22)$&$(\delta-1)(\delta-8)$&$-12$\\
\bottomrule
\end{tabular}
\end{center}

To calculate the scalars in the table, we use the following model. In each pair
fiber, the six coordinate rays of the perpendicular $\mathscr{D}_6$ system form a canonical six-element set.
Since the normalizer $N_W(X_0)$ preserves this six-element set and the CFW subgroup $A_{X_0}$ fixes each ray
individually, the induced permutation transport is independent of the choice of coset representatives.
This gives a decorated permutation model of dimension
$22680$, on which we check all $90720$ transport bridge choices, all $158$
defining relations, and all $960$ simple-reflection/root-contraction
conjugation identities. Applying the Specht modules of $S_6$
to these permutations, we obtain the integral pair-layer operators
\[
 T_{2,\theta}(\delta)=\delta^2I+\delta B_{\theta,1}+B_{\theta,0}.
\]
We evaluate the resulting scalar quadratics $q(\delta)$ both by stabilizer character-projection
traces and by central Lagrange interpolation on the same fiber model.

For the central Lagrange calculation, put $Z_1=\sum_{a>0}r_a$
and let $I_{X_0}$ be the base-fiber embedding. Let $\Sigma$ be the set
of distinct reflection-sum scalars given by the branching calculation.
We verify over $\mathbb Z$ that
$\prod_{c\in\Sigma}(Z_1-cI)I_{X_0}=0$, using normalizer-equivariance
to reduce the $3780$ fibers to $12$ orbits.
Since $I_{X_0}$ generates the induced module under $W(\mathscr{E}_8)$ and $Z_1$
is central, this polynomial annihilates the whole module.
We then form the integer matrix
\[
 V=g^{-1}\prod_{c\in\Sigma\setminus\{c_0\}}(Z_1-cI)I_{X_0}.
\]
Here $g$ is the positive greatest common divisor of the entries of
the unnormalized product, in the chosen integral fiber bases.
For the four choices $(\theta,c_0;a,b)$ given by
\[
 ([4,2],18;2,-8),\quad([2,2,1,1],-18;-4,-32),\quad
 ([3,2,1],8;7,-8),\quad([3,2,1],-8;-9,8),
\]
we verify over $\mathbb Z$ across all $3780$ fibers that
$Z_1V=c_0V$, $B_{\theta,1}V=aV$, and $B_{\theta,0}V=bV$, where the base blocks are $128I_9$, $64I_9$, $-528I_{16}$, and $-528I_{16}$, respectively.
Since its annihilating polynomial over $\mathbb Q$ is squarefree, $Z_1$ is
diagonalizable with spectrum contained in $\Sigma$, and
$\prod_{c\in\Sigma\setminus\{c_0\}}(Z_1-cI)$ is a nonzero scalar
multiple of the projection onto the $c_0$-eigenspace. Therefore, the $W(\mathscr{E}_8)$-translates
of the columns of $V$ span the $c_0$-eigenspace. Since
$B_{\theta,1}$ and $B_{\theta,0}$ commute with the $W(\mathscr{E}_8)$-action, we have
$B_{\theta,1}=aI$ and $B_{\theta,0}=bI$ on the entire $c_0$-eigenspace, where
$T_{2,\theta}(\delta)$ acts as the scalar $q(\delta)=\delta^2+a\delta+b$.
Since all positive-cardinality composition factors of the kernel $K_{2,\theta}$ have already been ruled out at $\delta_0\in\{-3,5\}$, every root contraction $e_a$
annihilates $K_{2,\theta}$, so $Z_1$ acts on any $s=0$ composition factor of $K_{2,\theta}$ by the corresponding eigenvalue $c_0$. Since $q(\delta_0)\ne0$ in all six cases, we get $K_{2,\theta}=0$.
This finishes the proof of Theorem~\ref{thm:app-exceptional}.

\subsection{Folded type \texorpdfstring{$\mathscr{F}_4$}{F4}}
\label{subsec:type-F4}
\label{subsec:f4-certified}

For the folded types, the ambient simply laced product still determines
the pairing, but several folded Weyl orbits may belong to one layer.

\begin{definition}\label{def:f4-model}
The \textbf{folded $\mathscr{F}_4$ algebra} is the subalgebra of
$\Br(\mathscr{E}_6;\delta)$ generated by the following elements, where the
ambient diagram has edges $13,34,45,56,24$:
\[
 \begin{aligned}
 r_1&=R_1R_6,& r_2&=R_3R_5,& r_3&=R_4,& r_4&=R_2,\\
 e_1&=E_1E_6,& e_2&=E_3E_5,& e_3&=E_4,& e_4&=E_2.
 \end{aligned}
\]
This is the realization in \cite[Theorem~1.2]{LiuF42013},
identified there with the presentation of Definition~1.1 of that paper.
\end{definition}

Here the short-root contractions $e_1,e_2$ satisfy $e_i^2=\delta^2e_i$, and the
long-root contractions $e_3,e_4$ satisfy $e_i^2=\delta e_i$. Under the diagram involution $\sigma=(1,6)(3,5)$, the twelve $\sigma$-fixed positive
roots of $\mathscr{E}_6$ form the long-root $\mathscr{D}_4$
subsystem of $\mathscr{F}_4$.

Liu's fixed-monomial basis and cellular stratification
\cite[Corollaries 5.1--5.2 and Theorem 6.5]{LiuF42013} have the following
layer data:
\[
\begin{array}{c|rrrr}
\text{cardinality of the admissible }\mathscr{E}_6\text{ set}&0&1&2&4\\
\text{number of fixed sets}&1&12&30&39\\
\text{coefficient group}&W(\mathscr{F}_4)&W(\mathscr{B}_3)&\mathfrak S_3&1.
\end{array}
\]
The pair layer has two Weyl-group orbits of sizes $12$ and $18$, and the
four-set layer has two orbits of sizes $36$ and $3$. These orbits belong
to one cellular layer in each case; they must not be separated when taking
the Gram determinant. In particular,
\[
 \dim\Br(\mathscr{F}_4;\delta)=1152+12^2\cdot48+30^2\cdot6+39^2=14985.
\]

Let $S_F$ be the union of the eight integers and the roots of the
five polynomials in the following statement. These polynomials are
pairwise coprime and squarefree, have no integer root, and have total
degree twelve. Hence $S_F$ has twenty elements.

\begin{theorem}\label{thm:f4-certified}
For $\delta\in\mathbb C^\times$, the folded algebra $\Br(\mathscr{F}_4;\delta)$ is
semisimple if and only if $\delta$ is neither one of
\[
 -8,-4,-2,-1,1,2,3,4
\]
nor a zero of any of the five polynomials
\begin{equation}\label{eq:f4-nonlinear}
\begin{split}
 &\delta^2+5\delta-8,\qquad \delta^2-3\delta+8,\qquad
       \delta^2-\delta-8,\\
 &\delta^3+7\delta^2+14\delta+32,\qquad
       \delta^3-2\delta^2+5\delta-16.
\end{split}
\end{equation}
Here $\Br(\mathscr{F}_4;\delta)$ is the complex specialization at
$z=\delta\in\mathbb C^\times$ of Liu's Laurent-parameter algebra over
$R_\times$. This statement makes no assertion about a polynomial form at zero.
\end{theorem}

\begin{proof}
By \cite[Corollaries~5.1--5.2 and Theorem~6.5]{LiuF42013}, the normal-form basis of $\Br(\mathscr{F}_4;\delta)$ consists
of the $\sigma$-fixed CFW normal forms of $\Br(\mathscr{E}_6;\delta)$, in which both admissible boundary sets are
$\sigma$-fixed and the middle element lies in the $\sigma$-fixed coefficient subgroup.
Since multiplication in $\Br(\mathscr{F}_4;\delta)$ is induced from $\Br(\mathscr{E}_6;\delta)$, we obtain the bilinear pairing in each folded layer by restricting the ambient CFW pairing to these $\sigma$-fixed normal forms and evaluating on the irreducible representations of the folded coefficient group. Passing between normalized idempotents and unnormalized root products
only multiplies the Gram determinants by nonzero constants and powers of $\delta\ne0$.

For the singleton layer ($t=1$), the twelve $\sigma$-fixed roots form the long-root $\mathscr{D}_4$
subsystem, and the triality permutations of the orthogonal roots give the semidirect factor
$\mathfrak S_3$, so the coefficient group is
$(C_2)^3\rtimes\mathfrak S_3=W(\mathscr{B}_3)$.
Under the action of $\mathfrak S_3$, the eight characters of $(C_2)^3$ form four orbits with representatives
$+++,-++,--+,---$. On the $12$-element $\mathscr{D}_4$ singleton layer, the four corresponding pairing matrices
have determinants
\begin{align*}
 f_0&=\delta^9(\delta-4)^2(\delta+8),&
 f_1&=(\delta-2)^8(\delta+4)^4,\\
 f_2&=\delta^6(\delta-4)^3(\delta+4)^3,&
 f_3&=(\delta-4)^4(\delta+2)^8.
\end{align*}
Here the characters $+++$ and $---$ have stabilizer $\mathfrak S_3$, whereas $-++$ and $--+$
have stabilizer $C_2$ and orbits of size three.
Restricting each of the ten irreducible representations of $W(\mathscr{B}_3)=(C_2)^3\rtimes\mathfrak S_3$ to $(C_2)^3$, we see that its $(C_2)^3$-isotypic summands are orthogonal
under the singleton pairing and permuted by triality. Therefore, up to nonzero scalar factors, the ten singleton Gram determinants of $\Br(\mathscr{F}_4;\delta)$ are
$f_0^{\dim\eta}$ and $f_3^{\dim\eta}$ for
$\eta\in\operatorname{Irr}(\mathfrak S_3)$, together with two copies
of $f_1^3$ and two copies of $f_2^3$.

For the pair layer, we consider the $30$ $\sigma$-fixed orthogonal pairs. At a base
pair $B_0$, the three orthogonal pairs in $B_0^\perp$ give a faithful
permutation realization of the coefficient group $\mathfrak S_3$.
We transport this ordered triple along the CFW root action: if a positive
root $a$ belongs to the current pair, then $E_a$ acts by $\delta$;
if $a$ is orthogonal to the pair, then $E_a$ moves it to a higher layer and acts as
zero in the pair quotient; otherwise, we choose $b$ in the pair with
$(a,b)\ne0$ and transport both the pair and its coefficient frame by
$r_b r_a$. To apply a folded generator, we apply its
commuting ambient factors successively; although the intermediate pair need not be
$\sigma$-fixed, the final nonzero pair is $\sigma$-fixed, so the $30$-pair space is
invariant under $\Br(\mathscr{F}_4;\delta)$. We construct the action on the $39$
fixed four-sets with trivial coefficient group in the same way.

By direct calculation with these integer polynomial matrices, we obtain
\begin{align}
\det G_{2,\mathbf1}
 &=\delta^{27}(\delta-2)^{15}(\delta+1)^9(\delta+2)^2
       (\delta^2+5\delta-8)^2
       (\delta^3+7\delta^2+14\delta+32),\label{eq:f4-trivial-det}\\
\det G_{2,\mathrm{sgn}}
 &=\delta^{27}(\delta-4)^2(\delta-3)^9(\delta+1)^2
       (\delta+2)^{18}(\delta^2-3\delta+8),\label{eq:f4-sign-det}\\
\det G_{2,\mathrm{std}}/3^{30}
 &=\delta^{54}(\delta-3)^9(\delta-1)^{22}(\delta+1)^9(\delta+4)^{12}
       (\delta^2-\delta-8)^4
       (\delta^3-2\delta^2+5\delta-16)^2,\label{eq:f4-standard-det}\\
\det G_4
 &=\delta^{72}(\delta-2)^{15}(\delta-1)^{38}
       (\delta+2)^{28}(\delta+4)^3.\label{eq:f4-four-det}
\end{align}
Here we equip the two-dimensional standard representation of $\mathfrak S_3$ with the invariant Gram matrix
$\left(\begin{smallmatrix}2&1\\1&2\end{smallmatrix}\right)$ of determinant $3$, which gives the factor $3^{30}$ in \eqref{eq:f4-standard-det} (and similarly gives $3^{12}f_0^2$ and
$3^{12}f_3^2$ for the corresponding singleton forms). We verify all fourteen determinant formulas
by $650$ exact integer evaluations via
Lemma~\ref{lem:exact-matrix-certificates}(3), using $61,61,121,$ and $157$ integer points for the four higher-layer
forms, respectively.
We also check all $868$ defining-relation and
contravariance matrix identities of degree at most six by evaluating them at the seven integers $\delta=-2,-1,0,1,2,3,4$.

Since the sum of the squared cell dimensions is $14985=\dim\Br(\mathscr{F}_4;\delta)$, these cell modules exhaust the cellular basis. By \eqref{eq:cellular-test}, taking the
union of the nonzero roots of the above determinants gives our result.
\end{proof}

\begin{remark}
The nonlinear factors in \eqref{eq:f4-nonlinear} are compatible with the
folding. For example, on the orbit-constant part of the trivial pair module,
the row-sum quotient of $G_{2,\mathbf1}$ is
\[
 \begin{pmatrix}
 \delta^2+3\delta+8&3\delta\\
 2\delta&\delta^2+4\delta
 \end{pmatrix},
\]
whose determinant is
$\delta(\delta^3+7\delta^2+14\delta+32)$.
Here this $2\times2$ matrix arises directly from the action on the orbit-constant subspace of the cell module $\Delta(2,\mathbf1)$.
\end{remark}

\subsection{Triality-folded type \texorpdfstring{$\mathscr{G}_2$}{G2}}
\label{subsec:G2-certified}
The triality-folded Brauer algebra of type $\mathscr{G}_2$ has two
distinct contraction exponents.

\begin{definition}\label{def:g2-model}
With node $2$ central in $\mathscr{D}_4$, the \textbf{triality-folded $\mathscr{G}_2$ algebra}
is the subalgebra of $\Br(\mathscr{D}_4;\delta)$ generated by
$r_1=R_2$, $r_2=R_1R_3R_4$, $e_1=E_2$, and $e_2=E_1E_3E_4$.
This is \cite[Theorem~6.1]{LiuG22015}, with the ambient nodes
$2,3$ interchanged and the folded indices relabelled as in
Definition~\ref{def:weyl-brauer}.
\end{definition}

It follows that $e_1^2=\delta e_1$ and $e_2^2=\delta^3e_2$.
Liu's normal forms and cellular stratification
\cite[Lemma~9.7 and Theorems~6.1, 10.1]{LiuG22015} give a group
quotient of dimension $12$, a layer of dimension $3^2\cdot2$, and a
minimal ideal of dimension $3^2$. The dimension is therefore $39$.
This is different from the $66$-dimensional admissible-partition algebra
of type $\mathscr{I}_2(6)$; see also \cite[Remark~10.2]{LiuG22015}.

\begin{proposition}\label{prop:G2-criterion}
For $\delta\ne0$, the triality-folded algebra
$\operatorname{Br}(\mathscr{G}_2;\delta)$ is semisimple exactly when
$\delta\notin\{-2,1\}$. \end{proposition}
\begin{proof}
In the middle layer, the three folded singleton boundaries have coefficient group $C_2$, whose two characters $\epsilon=\pm1$ give two
three-dimensional cell modules. Restricting the ambient $\mathscr{D}_4$ singleton action to the three roots in the folded orbit, we obtain the matrix realizations
\[
 r_1=\begin{pmatrix}1&0&0\\0&0&1\\0&1&0\end{pmatrix},\quad
 r_2=\epsilon\begin{pmatrix}0&1&0\\1&0&0\\0&0&1\end{pmatrix},\quad
 e_1=\begin{pmatrix}\delta&1&1\\0&0&0\\0&0&0\end{pmatrix},\quad e_2=0,
\]
where the central involution $(r_1r_2)^3$ acts as $\epsilon$. In these bases, the CFW normal-form pairing on each of the two middle-layer modules has Gram matrix
\[
 G=\begin{pmatrix}\delta&1&1\\1&\delta&1\\1&1&\delta\end{pmatrix}
\]
(if we instead use coset bases with a negative-character representative, we get the congruent matrix $\operatorname{diag}(1,-1,-1)G\operatorname{diag}(1,-1,-1)$).
On the three folded maximal frames, the coefficient group is trivial, and the ambient root action gives
\[
 r_1=\begin{pmatrix}0&1&0\\1&0&0\\0&0&1\end{pmatrix},\quad
 r_2=\begin{pmatrix}1&0&0\\0&0&1\\0&1&0\end{pmatrix},\quad
 e_1=\begin{pmatrix}0&0&0\\0&0&0\\1&1&\delta\end{pmatrix},\quad
 e_2=\begin{pmatrix}\delta^3&\delta^2&\delta^2\\0&0&0\\0&0&0\end{pmatrix}.
\]
For the three coset representatives $u,v\in\{1,r_1,r_2r_1\}$ in Liu's normal forms, each product $e_2u^{-1}ve_2$ is the corresponding pairing coefficient times $e_2$. Direct multiplication of the displayed matrices gives the minimal-layer form $\delta^2G$. Thus the entries of $G$ and $\delta^2G$ are the exact middle coefficients of the cellular basis products in Liu's normal forms.

Since $G$ has eigenvalues $\delta+2$ (with multiplicity $1$) and $\delta-1$ (with multiplicity $2$),
the three non-group cellular forms are nonsingular for $\delta\ne0$ if and only if
$\delta\notin\{-2,1\}$. The group-algebra cell forms of $\mathbb C[W(\mathscr{G}_2)]$ are nonsingular over
$\mathbb C$, and the sum of the squared cell dimensions is
$12+3\cdot3^2=39=\dim\Br(\mathscr{G}_2;\delta)$. By \eqref{eq:cellular-test}, we get our result.
\end{proof}

\subsection{Folded type \texorpdfstring{$\mathscr{H}_3$}{H3}}
\label{subsec:type-H3}
\label{app:h3}
The noncrystallographic types require their specified folded
presentations. We first recall the presentation and parameter
normalization for type $\mathscr{H}_3$.

\begin{definition}\label{def:h3-model}
\cite[Definition~2.1]{Liu2016}
Put $t=\delta^2$ with $\delta\ne0$. The \textbf{folded $\mathscr{H}_3$ algebra}
$\Br(\mathscr{H}_3;\delta)$ is generated by
$r_1,r_2,r_3,e_1,e_2,e_3$ with Coxeter diagram
$1\stackrel{5}{\text{---}}2\text{---}3$.
At each node $i$, we impose
\[
 r_i^2=1,\qquad r_ie_i=e_ir_i=e_i,\qquad e_i^2=t e_i.
\]
For the nonadjacent ordered pairs $(i,j)=(1,3),(3,1)$, we impose the
commutation relations between $r_i,r_j$, between $e_i,r_j$, and between $e_i,e_j$.
For the single-edge ordered pairs $(i,j)=(2,3),(3,2)$, we impose
\[
 r_ir_jr_i=r_jr_ir_j,\qquad
 r_jr_ie_j=e_ie_j,\qquad r_ie_jr_i=r_je_ir_j.
\]
At the five-bond edge, we impose
\begin{align*}
 r_1r_2r_1r_2r_1&=r_2r_1r_2r_1r_2,\\
 r_1r_2e_1r_2r_1&=r_2r_1e_2r_1r_2,\\
 e_1e_2e_1&=e_1,\qquad
 e_1r_2e_1=e_1,\qquad
 e_1r_2r_1r_2e_1=e_1.
\end{align*}
\end{definition}

This agrees with Chen's one-parameter presentation at
$\tau=t=\delta^2$, $\mu=1$
\cite[Definition~8.2 and Theorem~8.4]{Chen2012}; see also
\cite[Remark~5.4]{Liu2016}. The embedding into $\mathscr{D}_6$ and basis theorem
give dimension $1045$ \cite[Theorem~1.1]{Liu2016}.

\begin{theorem}\label{thm:app-h3}
For $\delta\in\C^*$, the folded algebra $\Br(\mathscr{H}_3;\delta)$ is
semisimple if and only if
\[
 \delta^2\notin\{-12,-4,-1,1,3,4\}.
\]
\end{theorem}
\begin{proof}
Write $W(\mathscr{H}_3)=A_5\times\langle c\rangle$, with $c$ central of order
two. Its fifteen reflections are the pairs $(a,c)$ with $a$ a double
transposition in $A_5$. Fix commuting reflections $s_0,s_2$ and put
$H=C_W(s_0)=\langle s_0,s_2,c\rangle\cong C_2^3$.
For $x,y\in\{1,-1\}$ let
\[
 \sigma_{x,y}(s_0)=1,\qquad
 \sigma_{x,y}(s_2)=x,\qquad \sigma_{x,y}(c)=y.
\]
In \cite[Section~7]{Chen2012}, Chen constructed four $15$-dimensional modules whose underlying $W(\mathscr{H}_3)$-representations are
$V_{x,y}=\Ind_H^W\sigma_{x,y}$, with basis $\{v_s\}$
indexed by the fifteen reflections. Their contraction operators are
\[
 E_sv_r=\begin{cases}
 t v_s,&r=s,\\
 0,&r\ne s,\ rs=sr,\\
 \rho(u)v_r,&rs\ne sr,
 \end{cases}
\]
where $u$ is the unique reflection satisfying $uru=s$ in the last case.
Write $E_sv_r=M_{x,y}(t)_{s,r}v_s$.
The fifth module has dimension five, indexed by the five commuting
reflection triples. Its distinct contraction row functionals form
$M_4(t)=(t-1)I_5+J_5$.
These module constructions are due to Chen \cite[Section~7]{Chen2012}.

Let $Z=\sum_{s\text{ reflection}}s$ be the sum of all fifteen reflections in $W(\mathscr{H}_3)$. From the action of $E_s$ above, we have
\[
 M_{x,y}(t)=tI+\rho(Z)-(1+x+xy)I.
\]
We use the following character values of $A_5$:
\[
\begin{array}{c|rrrrr}
\chi&\chi_1&\chi_3&\chi_{3'}&\chi_4&\chi_5\\\hline
\chi(1)&1&3&3&4&5\\
\chi(2A)&1&-1&-1&0&1.
\end{array}
\]
For $H\cap A_5\cong C_2^2$, the multiplicity of $\chi\in\Irr(A_5)$ in the induced representation is
$(\chi(1)+3\chi(2A))/4$ for the trivial character of $H\cap A_5$ and
$(\chi(1)-\chi(2A))/4$ for each of its three nontrivial characters.
Here $\sigma_{x,y}|_{H\cap A_5}$ is trivial if and only if $(x,y)=(1,1)$,
and $Z$ acts on an irreducible constituent $\chi$ by the scalar
$15y\chi(2A)/\chi(1)$. Substituting these eigenvalues, we obtain the determinants of the five matrices:
\[
\begin{array}{c|l}
\text{Matrix}&\text{Determinant}\\\hline
M_{1,1}&t^{10}(t-3)^4(t+12)\\
M_{-1,1}&(t-4)^6(t+1)^4(t+4)^5\\
M_{1,-1}&(t-4)^5(t-1)^4(t+4)^6\\
M_{-1,-1}&(t-4)^5(t-1)^4(t+4)^6\\
M_4&(t-1)^4(t+4).
\end{array}
\]

If one of these five matrices is singular, then the joint kernel of the contraction operators $E_s$ on that module
is a nonzero submodule with no invariant complement, since any invariant complement must contain the image lines $\operatorname{im}E_s$ and hence span
the entire module; in this case, $\Br(\mathscr{H}_3;\delta)$ is not semisimple.
Conversely, if all five matrices are invertible, then for any nonzero vector in one of these modules, some contraction $E_s$ maps it to a nonzero vector in a distinguished basis line, and by $W(\mathscr{H}_3)$-transitivity we see that each of the five modules is simple. Moreover, the four $15$-dimensional modules are pairwise nonisomorphic
because $H$ acts on $\operatorname{im}E_{s_0}$ via the distinct
characters $\sigma_{x,y}$, and the five-dimensional module is nonisomorphic to them
by dimension. The group quotient $\mathbb C[W(\mathscr{H}_3)]$ contributes
$|W|=120$ to the sum of squared simple module dimensions, and the five additional simple modules contribute
\[
 4\cdot15^2+5^2=925.
\]
Since $120+925=1045=\dim\Br(\mathscr{H}_3;\delta)$, this exhausts the dimension of the algebra and proves that $\Br(\mathscr{H}_3;\delta)$ is
semisimple when all five matrices are invertible. Since $t=\delta^2\ne0$, the nonzero roots of the five determinants give our result.
\end{proof}
\subsection{Folded type \texorpdfstring{$\mathscr{H}_4$}{H4}}
\label{subsec:type-H4}
\label{subsec:H4-status}

In type $\mathscr{H}_4$, we use the ambient $\mathscr{E}_8$ realization to construct both the algebra
and the coefficient transport for its cellular Gram matrices.

\begin{definition}\label{def:h4-model}
Put $t=\delta^2\ne0$. The \textbf{folded $\mathscr{H}_4$ algebra}
$A=\operatorname{Br}(\mathscr{H}_4;\delta)$ is the subalgebra of
$\operatorname{Br}(\mathscr{E}_8;\delta)$ generated by the folded elements below.
Here we number the simple roots of $\mathscr{E}_8$ by $0,\ldots,7$ with edges
$02,13,23,34,45,56,67$, and we take the admissible partition
\[
\begin{gathered}
 P_1=\{1,4\},\qquad P_2=\{3,5\},\qquad
 P_3=\{2,6\},\qquad P_4=\{0,7\},\\
 r_{\kappa}=\prod_{a\in P_{\kappa}}R_a,\qquad
 e_{\kappa}=\prod_{a\in P_{\kappa}}E_a.
\end{gathered}
\]
This is \cite[Theorem~1.2]{Liu2016}; our ambient node indices are
one less than those in that theorem.
\end{definition}

In particular $e_{\kappa}^2=t e_{\kappa}$. We describe the layer boundaries using the corresponding root orbits. Write $\mathcal X$
for the orbit of $P_1$,
viewed as $60$ orthogonal pairs of $\mathscr{E}_8$ root lines, and $\mathcal F$ for the
$75$ mutually orthogonal four-element subsets of $\mathcal X$.
Each $F\in\mathcal F$ is an $\mathscr{E}_8$ $8$-frame.  The pairs in $\mathcal X$
partition the $120$ $\mathscr{E}_8$ root lines.  Each $X\in\mathcal X$ belongs to five
frames; denote this set by $D(X)$.

Liu's normal forms and cellular stratification
\cite[Lemmas~3.1, 3.3, 3.6, Theorems~4.4 and~5.3]{Liu2016} give
\[
 \dim A=|W(\mathscr{H}_4)|+60^2|A_5|+75^2=236025.
\]
The group quotient is $\mathbb C[W(\mathscr{H}_4)]$.  The other cell modules are
$\Delta(\chi)$, of dimension $60\chi(1)$ for
$\chi\in\operatorname{Irr}(A_5)=\{\chi_1,\chi_{3,+},\chi_{3,-},\chi_4,\chi_5\}$,
and a maximal-layer module $\Delta_{\max}$ of dimension $75$.
The stabilizer $N(X)$ has order $240$, and its action on $D(X)$ has image
$A_5$ and kernel
$\langle r_X,c_{\mathscr{H}_3(X)}\rangle\cong C_2^2$.
Here $r_X=\prod_{a\in X}R_a$ is the product of the ambient root
reflections, and $c_{\mathscr{H}_3(X)}$ is the central involution of the
perpendicular folded subgroup of type $\mathscr{H}_3$.

We construct the layer representations and identify them with the cell
modules of $A$. We fix an ordering of $D(P_1)$, transport it to each $D(X)$ along
chosen words in the generators $r_{\kappa}$, and write $\pi(g;X)$ for the induced permutation
from $D(X)$ to $D(gX)$; for every simple generator, this permutation is even.
For nonorthogonal pairs $X\ne Y$, we apply the ambient $\mathscr{E}_8$ root-set rule
\[
 E_aB=R_bR_aB\qquad(b\in B,\ (a,b)\ne0)
\]
successively for the two roots in $Y$. We denote the resulting transport permutation from $D(X)$
to $D(Y)$ by $\pi_{Y,X}\in A_5$, which is independent of the
permitted choices of $b\in B$. We set
\begin{equation}\label{eq:H4-group-ring-form}
 K_{Y,X}(t)=
 \begin{cases}
 t\,1,&Y=X,\\
 0,&Y\perp X,\\
 \pi_{Y,X},&\text{otherwise}.
 \end{cases}
\end{equation}
On $\bigoplus_{X\in\mathcal X}\chi$, the generator $r_{\kappa}$ carries the
$X$-summand to the $r_{\kappa}X$-summand with coefficient
$\rho_\chi(\pi(r_{\kappa};X))$, whereas $e_{\kappa}$ has only the $P_{\kappa}$ row, with
coefficient $\rho_\chi(K_{P_{\kappa},X}(t))$.
For a nonsingular invariant form $J_\chi$ on $\chi$, put
$\mathcal J_\chi=\operatorname{diag}_{X\in\mathcal X}(J_\chi)$.
Then the cell form is
\begin{equation}\label{eq:H4-actual-pair-form}
 G_\chi(t)=\mathcal J_\chi\rho_\chi(K(t)).
\end{equation}

For the maximal layer, we set $e_F=e_Ue_V$, where $U$ and $V$ are any two distinct
root-pairs in $F$; by Liu's closure relation, $e_F$ is independent of the
choice of $U,V$ and satisfies $e_F^2=t^2e_F$. On the $75$ $\mathscr{E}_8$ $8$-frames, each ambient contraction $E_a$ sends a frame $B$
to $\delta B$ if $a\in B$ and to $R_bR_aB$ (for $b\in B$ with $(a,b)\ne0$) otherwise. From this action, we obtain both the $75$-dimensional representation of $A$
and the scalar coefficient of $e_F$ on each frame.

\begin{theorem}\label{thm:H4-complete-criterion}
The five pair-layer forms have the following monic determinants in $t$:
\begin{align*}
 \det_{\mathrm{mon}}G_{\chi_1}(t)
   &=(t-4)^{25}(t-1)^{16}(t+4)^{18}(t+44),\\
 \det_{\mathrm{mon}}G_{\chi_{3,+}}(t)
 =\det_{\mathrm{mon}}G_{\chi_{3,-}}(t)
   &=(t-16)^9(t-11)^{16}(t-8)^{25}(t+4)^{130},\\
 \det_{\mathrm{mon}}G_{\chi_4}(t)
   &=(t-13)^{25}(t-1)^{174}(t+11)^{25}(t+14)^{16},\\
 \det_{\mathrm{mon}}G_{\chi_5}(t)
   &=(t-4)^{216}(t+8)^{50}(t+11)^{16}(t+16)^{18}.
\end{align*}
More precisely, with
$Z=\sum_{u\in\operatorname{Ref}(W(\mathscr{H}_4))}u$ and
$V_\chi=\operatorname{Ind}_{N(P_1)}^{W(\mathscr{H}_4)}\chi$, where $\chi$ is inflated
through $N(P_1)\twoheadrightarrow A_5$, one has
\begin{equation}\label{eq:H4-casimir}
 G_\chi(t)=\mathcal J_\chi
 \bigl((t-1-c_\chi)I+\rho_{V_\chi}(Z)\bigr),
 \qquad c_\chi=\frac{15\chi(2A)}{\chi(1)}
       \in\{15,-5,-5,0,3\}.
\end{equation}

The maximal-layer form is
\begin{equation}\label{eq:H4-G75}
 G_{75}(t)=t^2I+tA_*+(J-I-A_*)
          =(t-1)((t+1)I+A_*)+J,
\end{equation}
where $A_*(F,F')=1$ exactly when the coefficient of $e_F$ on $F'$ is $t$.
This is a symmetric valency-$18$ graph, with spectrum
\[
 (-5)^{25},\quad(-2)^8,\quad3^{41},\quad18^1.
\]
Therefore, we have
\begin{equation}\label{eq:H4-det-75}
 \det G_{75}(t)=(t-1)^{82}(t-4)^{25}(t+4)^{42}(t+14)
\end{equation}
(if we use the unnormalized product of all eight ambient root contractions without dividing by $t^2$, the determinant is multiplied by $t^{150}$).
In particular, $A=\operatorname{Br}(\mathscr{H}_4;\delta)$ is semisimple if and only if
\begin{equation}\label{eq:H4-bad-set}
 t\notin\{-44,-16,-14,-11,-8,-4,1,4,8,11,13,16\}.
\end{equation}
\end{theorem}
\begin{proof}
We first verify that the above constructions realize the cell modules and cell forms of $A$.
By direct calculation with the root action, we generate the $60$-element orbit $\mathcal X$ from the folded reflections and construct
the $75$-element frame set $\mathcal F$ by enumerating mutually orthogonal quadruples in $\mathcal X$. We verify that the
signed-root group has order $14400$ and that the stabilizer $N(X)$ has order $240$, with image $A_5$ and kernel $C_2^2$ on $D(X)$.
We also check that all $3360$ permitted ambient transport choices on
nonorthogonal ordered pairs $(X,Y)$ yield the same element $\pi_{Y,X}\in A_5$ in
\eqref{eq:H4-group-ring-form}. Next, we verify all $42$ defining relations of $A$ as
monomial identities (with formal powers of $t$ and coefficients in $A_5$) on each of the $60$ pair columns and $75$ frame columns,
and we check the contravariance of the forms directly
in the group ring for all eight generators across all $28800$ pair entries and
$45000$ maximal-layer entries.

Let $I_{\max}$ be the two-sided ideal generated by $e_1e_3$. By Liu's normal forms, the idempotent
$\bar f=t^{-1}e_1$ in $\bar A=A/I_{\max}$ satisfies
$\bar f\bar A\bar f\cong\mathbb C[A_5]$, and $\bar A\bar f$ is free of
right rank $60$ over $\bar f\bar A\bar f$. Since our constructed pair module is generated
by its $P_1$ fiber (on which $\bar f\bar A\bar f$ acts by $\chi$) and is annihilated by $I_{\max}$, the induced homomorphism
\[
 \bar A\bar f\otimes_{\mathbb C[A_5]}\chi\longrightarrow
 \bigoplus_{X\in\mathcal X}\chi
\]
is surjective, and thus is an isomorphism by comparing dimensions. This identifies
our pair module with Liu's cell module $\Delta(\chi)$. Similarly,
$f_{\max}=t^{-2}e_1e_3$ has a one-dimensional corner $f_{\max}Af_{\max}\cong\mathbb C$ and $Af_{\max}$ has dimension
$75$, so the induced map onto our $75$-dimensional frame module is also an isomorphism. The
self-adjoint idempotents $t^{-1}e_X$ determine the cell form on $\Delta(\chi)$ from its
restriction $tJ_\chi$ to a single fiber, which gives
\eqref{eq:H4-actual-pair-form}. In the same way, $t^{-2}e_F$ determines the maximal-layer cell form $G_{75}(t)$, and we verify that all $33750$ choices of
two root-pairs $U,V\subset F$ give the same matrix entries.

For each $X\in\mathcal X$, the reflection $r_X$ acts trivially on $D(X)$, while the $15$
reflections orthogonal to $X$ act on $D(X)$ as the $15$ double transpositions of $A_5$.
The remaining $44$ reflections map $X$ bijectively onto its $44$ nonorthogonal
neighbors $Y$, and we check across all $2640$ ordered pairs that their induced permutations on the frames agree with $\pi_{Y,X}$ in
\eqref{eq:H4-group-ring-form}.
Subtracting the diagonal contributions from the reflection class sum $Z$, we obtain
\eqref{eq:H4-casimir}. Since $Z$ is a central class sum in
$\mathbb C[W(\mathscr{H}_4)]$, it acts on each irreducible constituent $\lambda\in\Irr(W(\mathscr{H}_4))$ by the scalar
$60\chi_\lambda(\mathrm{ref})/\chi_\lambda(1)$, which is both rational and an algebraic integer, hence a rational integer.

To compute the spectra of the pair forms, we realize the nontrivial irreducible representations of $A_5$
via its $4$-dimensional standard representation $\chi_4$, the
$6$-dimensional exterior square $\Lambda^2\chi_4\cong\chi_{3,+}\oplus\chi_{3,-}$, and the $5$-dimensional representation $\chi_5$ obtained
from the conjugation action on the six Sylow $5$-subgroups.
For $M_\chi=-\mathcal J_\chi^{-1}G_\chi(0)$, we verify the squarefree
minimal polynomials with root sets
\[
 \{-44,-4,1,4\},\quad\{-4,8,11,16\},\quad
 \{-14,-11,1,13\},\quad\{-16,-11,-8,4\},
\]
respectively, and we compute the eigenvalue multiplicities from the exact Lagrange projector traces
\[
 \operatorname{Tr}\left(\prod_{s\ne r}
             \frac{M_\chi-sI}{r-s}\right).
\]
On the exterior square $\Lambda^2\chi_4$, we construct $C$ from a central
element of $\mathbb Q[A_5]$. We verify $C^2=5I$ and
$\operatorname{Tr}C=0$, so $C$ acts by opposite scalars $\sqrt5$ and
$-\sqrt5$ on $\chi_{3,+}$ and $\chi_{3,-}$. Its block-diagonal extension
to the $60$ fibers commutes with $M_{\Lambda^2\chi_4}$. We verify that
the trace of its product with each eigenspace projector is zero.
Hence the eigenvalue multiplicities split equally between the two constituents.
Before using fixed-width arithmetic, the verification script checks bounds on the integer matrix products and traces to exclude overflow. These calculations are therefore exact.
In the same way, we verify the squarefree annihilating polynomial
of the adjacency matrix $A_*$ in \eqref{eq:H4-G75}, with roots
$-5,-2,3,18$. The corresponding projector traces give eigenvalue
multiplicities $25,8,41,1$, respectively. Since $A_*\mathbf1=18\mathbf1$, substituting these eigenvalues into
\eqref{eq:H4-G75} gives the four eigenvalue factors
$(t-1)(t-4)$, $(t-1)^2$, $(t-1)(t+4)$, and $(t+4)(t+14)$ of $G_{75}(t)$,
which proves \eqref{eq:H4-det-75}. Finally, the group quotient and the six
non-group cell modules exhaust the dimension of $A$:
\[
 |W(\mathscr{H}_4)|+60^2(1^2+3^2+3^2+4^2+5^2)+75^2=236025.
\]
By \eqref{eq:cellular-test}, the union of their roots gives
\eqref{eq:H4-bad-set}.
\end{proof}

\subsection{Dihedral types \texorpdfstring{$\mathscr{I}_2(m)$}{I2(m)}}
\label{subsec:type-I2}
\label{app:folded-additional}

We consider the natural four-bond diagram fold $\mathscr{I}_2(4)$, the odd dihedral diagram
algebras $\mathcal D_m(\delta)$ ($m\ge3$ odd), and the even dihedral diagram algebras
$\mathcal E_{2k}(\delta)$ ($k\ge3$). We then compare these models with Chen's
dihedral algebra.

The four-bond case is the symmetric $\mathscr{C}_2$ diagram model. We treat it
first, then the odd and even families with larger bond labels.

\begin{definition}\label{def:i24-model}
The \textbf{four-bond diagram algebra} $D$ is the
$\mathbb C[\delta]$-subalgebra of the ordinary four-strand Brauer algebra generated by
\[
 r_0=R_2,\quad r_1=R_1R_3,\quad e_0=E_2,\quad e_1=E_1E_3.
\]
\end{definition}

\begin{proposition}\label{prop:I24-C2}
The algebra $D$ has as a free $\mathbb C[\delta]$-basis all $25$ Brauer diagrams fixed by simultaneous reversal of the top and bottom indices.  In particular it is the polynomial symmetric-diagram model of type $\mathscr{C}_2$; its generator squares are $e_0^2=\delta e_0$ and $e_1^2=\delta^2e_1$.  Its complex specialization is semisimple exactly when $\delta\notin\{0,1,-2\}$.
\end{proposition}
\begin{proof}
Reversal preserves concatenation, so the generated algebra is
contained in the span of reversal-invariant diagrams. There are eight
such permutation diagrams: they form the group generated by $r_0,r_1$,
which is dihedral of order eight. A rank-two invariant diagram has one
self-conjugate horizontal pair in each row, with two choices per row,
and one of two bijections between the free pairs. Its eight possibilities
are represented by the words
\[
\begin{gathered}
 e_0,\quad r_1e_0,\quad e_0r_1,\quad r_1e_0r_1,\\
 e_0r_1r_0,\quad r_1e_0r_1r_0,\quad
 r_0r_1e_0,\quad e_0r_1r_0r_1.
\end{gathered}
\]
There are three invariant perfect matchings of one row of four vertices.
Thus there are nine rank-zero invariant diagrams, represented by
\[
 x e_1 y,\qquad x,y\in\{1,r_0,e_0\}.
\]
Each displayed word concatenates without a closed middle component;
within each list the resulting boundary matchings and free bijections
are distinct. These $8+8+9=25$ diagrams therefore belong to $D$ over
$\mathbb C[\delta]$, and the ambient diagram basis proves linear independence.
The two contraction squares count one and two loops, respectively.
This identifies $D$ with the polynomial symmetric-diagram model of
type $\mathscr{C}_2$.

For $\delta\ne0$, we apply Theorem~\ref{thm:C-semisimple} at $n=2$. Since $Z(2)=\{0\}$, the nonzero nonsemisimple parameters form the union of the nonzero
roots of $Q_+(\delta)=(\delta^2+\delta-2)/2$ and
$Q_-(\delta)=(\delta^2-\delta)/2$, namely $\delta=-2$ and $\delta=1$. At $\delta=0$, the span of the nine rank-zero symmetric diagrams forms a nonzero two-sided ideal of $D(0)$ whose square is zero, because concatenating any two rank-zero diagrams produces at least one closed middle loop. Therefore, $D(0)$ is not semisimple. This proves the criterion for all $\delta\in\mathbb C$.
\end{proof}

\label{subsec:odd-dihedral-spectrum}

For the odd family, we use the algebra obtained from Liu's admissible
partition of $\mathscr{A}_{m-1}$, with $m\geq3$ odd; for $m\geq5$ this is the algebra
of \cite[Definition~1.2 and Theorem~1.6]{LiuI22014}.
The diagram description below also specifies a polynomial form at
$\delta=0$, without specializing a Laurent coefficient ring at zero.

\begin{definition}\label{def:odd-dihedral}
Let $m\ge3$ be odd. Number the $m$ vertices by $\mathbb Z/m\mathbb Z$.
For $a\in\mathbb Z/m\mathbb Z$,
let $M_a$ be the partial matching whose free vertex is $a$ and whose edges are
\[
 \{a+j,a-j\},\qquad 1\leq j\leq(m-1)/2.
\]
Write $D_{a,b}$ for the Brauer diagram with top matching $M_a$, bottom
matching $M_b$, and its unique through strand joining the free vertices.
The \textbf{odd dihedral diagram model} $\mathcal D_m(\delta)$
is the span of the $2m$ permutation diagrams of
the dihedral group and the $m^2$ diagrams $D_{a,b}$, inside the ordinary
Brauer diagram algebra on $m$ strands.
\end{definition}

As we verify in the proof of Theorem~\ref{thm:odd-dihedral-spectrum} below, these $2m+m^2$ diagrams form the basis
of a subalgebra of $B_m(\delta)$ for every $\delta\in\mathbb C$.
Up to relabelling the vertices by a Coxeter cycle, the two dihedral reflection
generators and their matching diagrams coincide with the images of Liu's
four folded generators, so $\mathcal D_m(\delta)$ is isomorphic to Liu's folded dihedral algebra
for $m\ge5$ and $\delta\ne0$ \cite[Theorem~1.6]{LiuI22014}.
When $m=3$, $\mathcal D_3(\delta)$ is the ordinary
three-strand Brauer algebra $B_3(\delta)$.

Since the layer pairing is invariant under cyclic translation, its Gram matrix is circulant. To express its
eigenvalues over $\mathbb Z[\delta]$, we group the Fourier
frequencies according to their conductors using Ramanujan sums.

\begin{definition}\label{def:ramanujan}
For a positive integer $q$, the \textbf{Ramanujan sum}
\cite[Sections~1--2, equations~(2.5)--(2.7)]{Ramanujan1918} is
\[
 c_q(j)=\sum_{\substack{0\leq a<q\\(a,q)=1}}
             \exp(2\pi\mathrm i aj/q)
       =\sum_{h\mid(q,j)}h\,\mu(q/h),
\]
where $\mu$ is the number-theoretic M\"obius function. For $q\mid m$, put
\begin{equation}\label{eq:dihedral-ramanujan}
 P_{m,q}(\delta)=\sum_{d\mid m}
       \delta^{(d-1)/2}c_{m/d}(m/q).
\end{equation}
We call $P_{m,q}$ the \textbf{odd-dihedral frequency polynomial}
for $q$.
\end{definition}

The divisor-sum formula follows by inserting
$[\gcd(a,q)=1]=\sum_{h\mid\gcd(a,q)}\mu(h)$ into the exponential
sum. In particular, these frequency polynomials lie in $\mathbb Z[\delta]$.

\begin{theorem}\label{thm:odd-dihedral-spectrum}
For odd $m\geq3$, the algebra $\mathcal D_m(\delta)$ has dimension $2m+m^2$.
Its non-group cellular form has matrix
\begin{equation}\label{eq:dihedral-gcd-gram}
 G_m(\delta)_{a,b}
       =\delta^{(\gcd(a-b,m)-1)/2},
       \qquad a,b\in\mathbb Z/m\mathbb Z,
\end{equation}
where $\gcd(0,m)=m$. In particular,
\begin{equation}\label{eq:dihedral-det}
 \det G_m(\delta)=\prod_{q\mid m}P_{m,q}(\delta)^{\varphi(q)},
 \qquad
 \rho_m(\delta)\coloneqq\operatorname{rank}G_m(\delta)
   =\sum_{\substack{q\mid m\\P_{m,q}(\delta)\ne0}}\varphi(q).
\end{equation}
Here $\varphi$ is Euler's totient function. Therefore,
$\mathcal D_m(\delta)$ is semisimple if and only if $P_{m,q}(\delta)\ne0$ for every $q\mid m$.

\end{theorem}

\begin{proof}
The two involutions on the middle vertices in the product
$D_{a,b}D_{c,d}$ are $x\mapsto2b-x$ and $x\mapsto2c-x$.
Their product is translation by $2(b-c)$. Put $g=\gcd(b-c,m)$.
Modulo $g$, the two involutions coincide, and their common action on residue
classes has the unique fixed class $b\equiv c\pmod g$.
On that class the overlay is the unique open path from $b$ to $c$; on each
of the $(g-1)/2$ other pairs of residue classes it is one closed component.
Indeed, on each residue class the translation has a single orbit of length
$m/g$, so none of these components splits further. The same description
when $b=c$ consists of the isolated free vertex and $(m-1)/2$ doubled edges.
It follows that
\[
 D_{a,b}D_{c,d}
   =\delta^{(\gcd(b-c,m)-1)/2}D_{a,d}.
\]
Since the dihedral permutations permute the matchings $M_a$, left and right multiplication by dihedral permutations
preserves the span of the $D_{a,b}$. Because the $2m+m^2$ basis diagrams are
distinct, this gives the dimension formula and the cellular multiplication rule.

The Gram matrix $G_m(\delta)$ in \eqref{eq:dihedral-gcd-gram} is circulant, and its Fourier eigenvector
of frequency $j\in\mathbb Z/m\mathbb Z$ has eigenvalue
\[
 \sum_{a\bmod m}\delta^{(\gcd(a,m)-1)/2}
                     \exp(2\pi\mathrm i aj/m)
   =\sum_{d\mid m}\delta^{(d-1)/2}c_{m/d}(j).
\]
This eigenvalue depends only on $q=m/\gcd(j,m)$ and equals $P_{m,q}(\delta)$. Since there are
$\varphi(q)$ frequencies $j$ with $\gcd(j,m)=m/q$ (including $j=0$ for
$q=1$), we obtain both formulas in \eqref{eq:dihedral-det}. We also note that $G_m(\delta)$ is never the zero matrix, since its entries with
$\gcd(a-b,m)=1$ are equal to $1$.

Let $I$ be the diagram ideal spanned by the $D_{a,b}$.
It is a cellular layer with Gram matrix $G_m(\delta)$, and
$\mathcal D_m(\delta)/I\cong\mathbb C[W(\mathscr{I}_2(m))]$.
The criterion follows from \eqref{eq:cellular-test}.
\end{proof}

\begin{corollary}\label{cor:odd-dihedral-zero}
The algebra $\mathcal D_m(0)$ is semisimple if and only if $m$ is
squarefree. For an odd prime $p$,
\[
 \det G_p(\delta)=
 (\delta^{(p-1)/2}-1)^{p-1}(\delta^{(p-1)/2}+p-1).
\]
\end{corollary}
\begin{proof}
At $\delta=0$, only the term $d=1$ survives in
\eqref{eq:dihedral-ramanujan}, so the classical Ramanujan-sum formula gives
\[
 P_{m,q}(0)=c_m(m/q)=\mu(q)\frac{\varphi(m)}{\varphi(q)}.
\]
Thus $P_{m,q}(0)\ne0$ for all $q\mid m$ if and only if $\mu(m)\ne0$, namely when $m$ is squarefree.
When $m=p$ is an odd prime, every off-diagonal entry of $G_p(\delta)$ is equal to $1$, which gives
the displayed determinant formula.
\end{proof}

\subsubsection*{Even dihedral diagram algebras: a polynomial form and its spectrum}

For even folds there are two non-group diagram ranks, so the Fourier
calculation involves three forms.
Liu's even dihedral algebra \cite[Definition~1.1]{LiuI22014}
is defined over $\mathbb Z[\delta^{\pm1}]$,
so its original coefficient ring has no specialization at $\delta=0$.
We use the following polynomial diagram form. Put $m=2k\ge6$.
In the ordinary Brauer algebra on $2k$ strands, let
\[
 r_0=\prod_{\substack{1\le \kappa<2k\\\kappa\ \mathrm{even}}}R_{\kappa},
 \quad r_1=\prod_{\substack{1\le \kappa<2k\\\kappa\ \mathrm{odd}}}R_{\kappa},
 \quad e_0=\prod_{\substack{1\le \kappa<2k\\\kappa\ \mathrm{even}}}E_{\kappa},
 \quad e_1=\prod_{\substack{1\le \kappa<2k\\\kappa\ \mathrm{odd}}}E_{\kappa},
 \qquad \rho=r_1r_0,\quad c=\rho^k.
\]
Here each product consists of commuting factors.  The generators $e_0,e_1$
have square parameters $\delta^{k-1},\delta^k$, respectively.  This is
the partition in \cite[Theorem~1.6]{LiuI22014}; it is distinct
from the $\mathscr{D}_4$ folding defining the $39$-dimensional algebra of type $\mathscr{G}_2$.

We relabel the vertices cyclically by $\mathbb Z/2k\mathbb Z$ so that left
multiplication by $\rho$ shifts the index of each top matching by $+1$.  For
$a\in\mathbb Z/k\mathbb Z$, let $M_a$ be the partial matching associated
to the involution $x\mapsto2a-x$, with free vertices $a,a+k$.  Let $F_a$
be the perfect matching associated to $x\mapsto2a+1-x$, and put
\[
 H_a=M_a\cup\{\{a,a+k\}\},\qquad
 \mathcal P=\{F_a:a\bmod k\}\sqcup\{H_a:a\bmod k\}.
\]
Thus $e_0$ has top and bottom $M_0$, $e_1$ has top and bottom $F_0$,
and $e_0e_1$ has top $H_0$ and bottom $F_0$, with no scalar factor.
For $k\ge3$ the $2k$ elements of $\mathcal P$ are distinct.  For odd $k$,
$H_a$ has a unique diameter, determining $a$, and its other edges have
sum $2a$; no $F_b$ with that diameter has this edge sum.  For even $k$,
$H_a$ has two diameters whereas $F_b$ has none; the nondiameter edges
distinguish $H_a$ from $H_{a+k/2}$, the only other candidate with those
same diameters (here $k\ge4$).

\begin{definition}\label{def:even-dihedral}
With the matchings just specified, the \textbf{even dihedral polynomial model}
$\mathcal E_{2k}$ is the $\mathbb Z[\delta]$-span of the following
ordinary Brauer diagrams:
the $4k$ dihedral permutation diagrams;
the $2k^2$ diagrams with top $M_a$, bottom $M_b$, and either bijection
between their two free vertices; and
the $4k^2$ diagrams with independently chosen top and bottom in $\mathcal P$.
Its complex fiber is
$\mathcal E_{2k}(z)=\mathcal E_{2k}\otimes_{\mathbb Z[\delta],\,\delta\mapsto z}\mathbb C$.
\end{definition}

These diagrams form a basis of a subalgebra, generated by the four displayed
generators.  Indeed, dihedral permutations preserve the three matching
families.  A product of two rank-two diagrams either preserves two through
strands or joins both free vertices on each boundary, producing an $H/H$
diagram.  A product involving a rank-zero diagram has rank zero; whenever
a boundary initially has matching $M_a$, its free vertices are joined,
giving $H_a$.  All products therefore remain in the indicated span.
Conversely all basis diagrams are represented, without any powers of
$\delta$ as scalar factors, by
\[
 \rho^a e_0c^\epsilon\rho^{-b},\qquad
 \rho^a v\rho^{-b}\quad
 (v=e_1,e_0e_1,e_1e_0,e_0e_1e_0),
 \quad 0\le a,b<k,\quad\epsilon=0,1.
\]
Since distinct ordinary Brauer diagrams are linearly independent,
$\mathcal E_{2k}$ is free of rank $4k+6k^2$ over $\mathbb Z[\delta]$, and by \cite[Theorem~1.6]{LiuI22014} its localization at $\delta$ is isomorphic to Liu's folded algebra. All results at $\delta=0$ for even dihedral types below refer
to $\mathcal E_{2k}(0)$.

For the spectral formulas, $[P]$ denotes $1$ if the assertion $P$ holds
and $0$ otherwise. The frequency $j$ modulo $k$ has order, also called
its conductor, $q=k/\gcd(j,k)$.

\begin{theorem}
\label{thm:even-dihedral-spectrum}
Write $k=Ln$, where $L=2^s$ and $n$ is odd.  The three non-group cell
modules of $\mathcal E_{2k}(\delta)$ have dimensions $k,k,2k$.
Put $g=\gcd(b-a,k)$, including $g=k$ when $a=b$.
Their first two Gram matrices satisfy
\[
 (G_{1,+})_{a,b}=
 \begin{cases}\delta^{g-1},&k/g\text{ odd},\\0,&k/g\text{ even},\end{cases}
 \qquad
 G_{1,-}=D G_{1,+}D,\qquad
 D_{a,a}=(-1)^{\lfloor a/L\rfloor},\quad 0\le a<k.
\]
In particular the matrices are integrally congruent.  Define
\[
 Q_{k,q}^{(1)}(\delta)
   =\sum_{d\mid n}\delta^{Ld-1}c_{n/d}(n/q)\qquad(q\mid n).
\]
Then
\begin{equation}\label{eq:even-dihedral-G1}
 \det G_{1,\pm}=\prod_{q\mid n}
       \bigl(Q_{k,q}^{(1)}\bigr)^{L\varphi(q)},\qquad
 \rho_1\coloneqq\operatorname{rank}G_{1,\pm}
       =L\sum_{\substack{q\mid n\\Q_{k,q}^{(1)}(\delta)\ne0}}\varphi(q).
\end{equation}

We order the basis of the third cell module first by $F_0,\ldots,F_{k-1}$ and then by $H_0,\ldots,H_{k-1}$,
with representatives $\rho^a e_1$ and $\rho^a e_0e_1$ ($0\le a<k$).
Then its Gram matrix is
\[
 G_0=\begin{pmatrix}A&B\\B^{\mathsf T}&C\end{pmatrix},\qquad
 A_{a,b}=\delta^g,\quad
 B_{a,b}=\delta^{(\gcd(2(b-a)-1,k)+1)/2},\quad
 C_{a,b}=\delta^{g+[k/g\ \mathrm{even}]}.
\]
For $q\mid k$, set
\[
 \alpha_{k,q}=\sum_{d\mid k}\delta^d c_{k/d}(k/q),\qquad
 \gamma_{k,q}=\sum_{d\mid k}\delta^{d+[k/d\ \mathrm{even}]}c_{k/d}(k/q),
\]
and for odd $q\mid n$ set
\[
 \beta_{k,q}=L\sum_{d\mid n}\delta^{(d+1)/2}c_{n/d}(n/q).
\]
Let
\[
 T_{k,q}=\begin{cases}
  \operatorname{diag}(\alpha_{k,q},\gamma_{k,q}),&q\text{ even},\\
  \begin{pmatrix}\alpha_{k,q}&\beta_{k,q}\\
                  \beta_{k,q}&\gamma_{k,q}\end{pmatrix},&q\text{ odd}.
 \end{cases}
\]
Under a constant invertible change of basis, $G_0$ is similar to the direct sum
of $\varphi(q)$ copies of $T_{k,q}$ over all $q\mid k$.  Hence
\begin{equation}\label{eq:even-dihedral-G0}
 \det G_0=\prod_{q\mid k}(\det T_{k,q})^{\varphi(q)},\qquad
 \rho_0\coloneqq\operatorname{rank}G_0
       =\sum_{q\mid k}\varphi(q)\operatorname{rank}T_{k,q}.
\end{equation}
For odd $q$, the block rank is two if its determinant is nonzero, zero
if $\alpha_{k,q}=\beta_{k,q}=\gamma_{k,q}=0$, and one otherwise.
For even $q$, it is the number of nonzero entries among
$\alpha_{k,q},\gamma_{k,q}$.
\end{theorem}

\begin{proof}
We give the loop counts explicitly.  In the overlay of $M_a$ and $M_b$,
the two reflection involutions have product a translation by $2(b-a)$.
There are $2g$ translation orbits, namely the residue classes modulo $2g$,
each of size $k/g$.  Reflection pairs all but the two classes
$a$ and $a+g$ modulo $2g$; each of the $g-1$ paired orbits gives one
closed alternating component.  Each exceptional class gives one path.
Write $t=(b-a)/g$ and $u=k/g$, so $\gcd(t,u)=1$ unless $a=b$,
whose doubled-edge overlay gives the same conclusions directly.
If $u$ is odd, each path has one endpoint free in $M_a$ and one free in
$M_b$.  The endpoint $a$ is connected to $b$ when $t$ is even and to
$b+k$ when $t$ is odd: these assertions follow by reducing the endpoints
modulo $2g$.  Thus the rank-two product has $g-1$ loops and central factor
$c^{t\bmod2}$.  If $u$ is even, $t$ is odd; the exceptional paths instead
join the two free vertices of $M_a$ to one another and those of $M_b$
to one another.  The product has rank zero and vanishes in the middle
cell quotient. This gives the formula for $G_{1,+}$. On its nonzero entries,
$g/L$ is odd, so
\[
 (-1)^{(b-a)/g}=(-1)^{(b-a)/L}
       =(-1)^{\lfloor a/L\rfloor+\lfloor b/L\rfloor}.
\]
Applying the two characters of $\langle c\rangle$, we obtain the congruence $G_{1,-}=DG_{1,+}D$.
Since the nonzero entries of $G_{1,+}$ occur only between indices in the same residue class modulo
$L$, $G_{1,+}$ decomposes into $L$ copies of the $n\times n$ circulant matrix with entries
$\delta^{L\gcd(b-a,n)-1}$. Diagonalizing this circulant by the same Ramanujan-sum calculation as in Theorem~\ref{thm:odd-dihedral-spectrum}, we get \eqref{eq:even-dihedral-G1}.

Next, in the overlay of $F_a$ and $F_b$, the product of the two reflection involutions is a translation
with $2g$ cycles, which pair up into $g$ closed alternating components. When we complete both
partial matchings $M_a$ and $M_b$ to $H_a$ and $H_b$ by adding the diameters $\{a,a+k\}$ and $\{b,b+k\}$, the two open paths in the overlay of $M_a$ and $M_b$ join
into a single loop if $u$ is odd, and close into two separate loops if $u$ is
even; this gives $g$ and $g+1$ loops for the overlay of $H_a$ and $H_b$, respectively.
For the mixed overlay of $F_a$ and $M_b$, the product of the reflection
involutions is translation by $2(b-a)-1$. Put
$h=\gcd(2(b-a)-1,2k)=\gcd(2(b-a)-1,k)$, which is odd.
Among the $h$ residue-class orbits, $h-1$ are paired by reflection into $(h-1)/2$ closed loops, and the remaining orbit forms a single open path joining the
free vertices $b$ and $b+k$. When we add the diameter $\{b,b+k\}$ of $H_b$, this open
path closes into one additional loop, which gives $(h+1)/2$ loops in total. Thus we obtain all three blocks $A,B,C$ of $G_0$.

To see that $G_{1,\pm}$ and $G_0$ are the cellular Gram matrices of $\mathcal E_{2k}(\delta)$, let $I_0$ be the ideal spanned by all
rank-zero diagrams and $I_1$ the ideal spanned by all non-permutation diagrams.
The basis of $I_0$ consists of the diagrams $C_{P,Q}$ for $P,Q\in\mathcal P$,
which multiply by
$C_{P,Q}C_{R,S}=(G_0)_{Q,R}C_{P,S}$.
Modulo $I_0$, the quotient $I_1/I_0$ has basis
$\rho^a e_0c^\epsilon\rho^{-b}$, which is an inflation of
$\mathbb C[C_2]$ whose two characters give $G_{1,+}$ and $G_{1,-}$.
Diagram reflection exchanges the two boundary indices and inverts
the middle group element. Together with a cellular basis of the dihedral
group algebra, these multiplication rules define a cellular datum for
$\mathcal E_{2k}(\delta)$ at every $\delta\in\mathbb C$.

Finally, all three $k\times k$ blocks $A,B,C$ are circulant. Using the Fourier basis vectors
$(\zeta_k^{-ja})_{a\bmod k}$, we see that $A$ and $C$ act on frequency $j$ by the scalars $\alpha_{k,q}$ and $\gamma_{k,q}$, where
$q=k/\gcd(j,k)$. Since the first row of $B$
is periodic modulo $n$, its eigenvalue on frequency $j$ vanishes unless $j=L\ell$ (that is, when $q$ is odd).
For $j=L\ell$, we set $v=(n+1)/2$; since $2v\equiv1\pmod n$ and $\gcd(2t-1,n)=\gcd(t-v,n)$, we have
\[
 \lambda_j(B)=\zeta_n^{-\ell v}\beta_{k,q},\qquad
 \lambda_j(B^{\mathsf T})=\zeta_n^{\ell v}\beta_{k,q}.
\]
Removing the reciprocal phase factors $\zeta_n^{\mp\ell v}$ in each $2\times2$ frequency block by a diagonal similarity transformation, we reduce $G_0$ on frequency $j$ to $T_{k,q}$.
Since there are $\varphi(q)$ frequencies of conductor $q$, we obtain
\eqref{eq:even-dihedral-G0} and the rank formulas. In terms of a primitive
$2k$th root of unity $\zeta_{2k}$, the phase factor is
$(-1)^\ell\zeta_{2k}^{-j}$, whose sign cancels in $\det T_{k,q}$.
\end{proof}

\begin{corollary}\label{cor:even-dihedral-structure}
For $k\ge3$ and every $\delta\in\mathbb C$,
$\mathcal E_{2k}(\delta)$ is semisimple if and only if
$Q^{(1)}_{k,q}(\delta)\ne0$ for all $q\mid n$ and
$\det T_{k,q}(\delta)\ne0$ for all $q\mid k$.
In particular, $\mathcal E_{2k}(0)$ is not semisimple.
\end{corollary}
\begin{proof}
The three non-group cell modules and the dihedral group quotient exhaust
the basis of $\mathcal E_{2k}(\delta)$, since $k^2+k^2+(2k)^2+4k=6k^2+4k$.
Applying \eqref{eq:cellular-test} to
\eqref{eq:even-dihedral-G1} and \eqref{eq:even-dihedral-G0}, we obtain the first assertion.
At $\delta=0$, we have $G_0(0)=0$ (or equivalently, the nonzero
rank-zero ideal $I_0$ satisfies $I_0^2=0$), so $\mathcal E_{2k}(0)$ is not semisimple.
\end{proof}

\subsubsection*{Comparison with Chen\textquotesingle s dihedral algebra}
\label{subsec:chen-dihedral-comparison}

The preceding criteria depend on the chosen folded model, not only on
the dihedral Coxeter group. To make the comparison precise, we record
the rank-two specialization of Chen's definition.

\begin{definition}\label{def:chen-dihedral}
Let $W=\mathscr{I}_2(m)$ with $m\ge3$, and let $\mathcal H$ be its reflecting
lines in its two-dimensional reflection representation. Choose loop
parameters $\tau_H$ constant on $W$-orbits of lines and nonzero
parameters $\mu_s$ constant on conjugacy classes of reflections.
\textbf{Chen's dihedral algebra} is generated by $\mathbb C[W]$ and
elements $e_H$ for $H\in\mathcal H$, subject to
\[
\begin{gathered}
 w e_H w^{-1}=e_{wH},\qquad e_H^2=\tau_H e_H,\\
 h e_H=e_H=e_Hh\quad(h\in N_W(H)),\\
 e_He_K=\left(\sum_{\substack{s\text{ a reflection}\\sK=H}}
                  \mu_s s\right)e_K\quad(H\ne K),
\end{gathered}
\]
where an empty sum is zero. Here $N_W(H)$ is the setwise stabilizer
of $H$. These are the dihedral relations of
\cite[Definition~1.1]{Chen2012}; our $\tau_H$ is Chen's parameter $m_H$.
\end{definition}

The normalizer absorption follows from relations (1) and (1)' of
\cite[Definition~1.1]{Chen2012}. For odd $m$, we have
$N_W(H)=\{1,s_H\}$. For even $m$, let $c=-I$ and $u=cs_H$.
The element $u$ is the reflection whose fixed line is perpendicular
to $H$. Thus $H\cap\operatorname{Fix}(u)=\{0\}$ is a noncrossing
edge, and relation (1)' gives $ue_H=e_H=e_Hu$. Together with
$s_He_H=e_H=e_Hs_H$, this also gives $ce_H=e_H=e_Hc$.
Since $N_W(H)=\{1,s_H,u,c\}$, the stated absorption follows.

Chen's cellular structures and dimensions are given in
\cite[Section~6, Theorems~6.1--6.3]{Chen2012}. For odd $m$, normalize its
reflection parameter to $\mu=1$ and write $\tau$ for its loop parameter.
The algebra has dimension $2m+m^2$, with one non-group cellular form
\[
 H_m(\tau)=(\tau-1)I_m+J_m,\qquad
 \det H_m(\tau)=(\tau-1)^{m-1}(\tau+m-1).
\]
Thus its bad parameters are exactly $1$ and $1-m$; in particular,
$\tau=0$ is semisimple. Chen's loop parameter must not be identified
with the folded parameter without comparing the defining relations.

For even $m=2k$, Chen's algebra has dimension $2m+m^2/2$ and exactly two
non-group cells, each of dimension $k$. It does not have two independent
sign choices for each of these cells. In the uniform-reflection
specialization $\mu_0=\mu_1=\mu$, their forms are
\[
 H_{k,\kappa}(\tau_{\kappa})=(\tau_{\kappa}-2\mu)I_k+2\mu J_k,
 \qquad \kappa=0,1,
\]
with determinants
\[
 (\tau_{\kappa}-2\mu)^{k-1}\bigl(\tau_{\kappa}+2\mu(k-1)\bigr).
\]
These formulas follow directly from Chen's multiplication: for distinct
hyperplanes in one reflection orbit, the coefficient is the sum of the
two reflection parameters. For $\mu=1/2$, the cell form indexed by
$\kappa$ is singular exactly when $\tau_{\kappa}\in\{1,1-k\}$.
More generally, Definition~\ref{def:chen-dihedral} allows two independent reflection
parameters in addition to the two loop parameters. In particular, the
$30$-dimensional Chen algebra for $m=6$ is distinct from the
$39$-dimensional $\mathscr{D}_4$-folded algebra and the $66$-dimensional
$\mathscr{A}_5$-folded algebra.

\section{Semisimplicity at zero and polynomial forms}
\label{sec:further-results}\label{subsec:master-zero}
Theorem~\ref{thm:master-nonzero-semisimplicity} gives the semisimplicity criteria for nonzero
parameters $\delta\in\mathbb C^\times$. We now consider the polynomial diagram forms defined above
at $\delta=0$ and complete their semisimplicity criteria over $\mathbb C$.
Since there is no unital ring homomorphism $\mathbb C[z,z^{-1}]\to\mathbb C$
sending $z$ to $0$, defining the zero fiber $A(0)$ of a Laurent-parameter Brauer algebra requires choosing a polynomial form in the sense
of Definition~\ref{def:polynomial-form}. This remark applies both to the
simply laced CFW presentation and to the folded Laurent presentations in
\cite{CFW2009,CohenLiu2014,LiuF42013,LiuG22015,Liu2016}.
In Proposition~\ref{thm:g2-dual-zero-forms} below, we construct two polynomial forms
of the triality-folded $\mathscr{G}_2$ Laurent algebra such that one is nonsemisimple at $\delta=0$ while the other is semisimple at $\delta=0$. Thus semisimplicity at zero depends on the polynomial form.

\subsection{The polynomial diagram families}
\begin{theorem}\label{thm:master-zero-semisimplicity}
For the polynomial diagram forms used in this paper, the answers at
zero are as follows.
\begin{enumerate}
\item $B_n(0)$, $n\ge1$, is semisimple exactly when $n\in\{1,3,5\}$.
\item The symmetric type $\mathscr{C}_n$ diagram algebra $A_n(0)$ is nonsemisimple
for every $n\ge1$.
\item For odd $m\ge3$, $\mathcal D_m(0)$ is semisimple exactly when
$m$ is squarefree.
\item The $\mathscr{I}_2(4)$ diagram fold and $\mathcal E_{2k}(0)$, $k\ge3$,
are nonsemisimple.
\end{enumerate}
\end{theorem}
\begin{proof}
The first assertion is the ordinary criterion
\cite{Rui2005,RuiSi2006}. For the second, let $I$ be the span of
symmetric diagrams with no through strands. It is a nonzero two-sided
ideal, because concatenation cannot increase through rank. The
concatenation of any two of its basis diagrams has a closed middle
component, so $I^2=0$ at zero. Hence $I\subseteq\operatorname{Rad}A_n(0)$.
The third assertion is Corollary~\ref{cor:odd-dihedral-zero}; the last
follows from Proposition~\ref{prop:I24-C2} and
Corollary~\ref{cor:even-dihedral-structure}.
\end{proof}

Combining Theorem~\ref{thm:master-zero-semisimplicity} with Theorem~\ref{thm:master-nonzero-semisimplicity}, we obtain the complete complex
semisimplicity criteria for these polynomial diagram families. Explicitly, for $n\ge2$, we set
\[
 \widehat Z(n)=(Z(n)\setminus\{0\})\cup
 \begin{cases}\varnothing,&n=3,5,\\\{0\},&n\ne3,5.\end{cases}
\]
Then $B_n(\delta)$ is semisimple exactly when
$\delta\notin\widehat Z(n)$, and $B_1(\delta)=\mathbb C$ for all
$\delta$. For $n\ge2$, the criterion
$Q_+(\delta),Q_-(\delta)\notin Z(n)$ remains valid for all
$\delta\in\mathbb C$, since $Q_-(0)=0\in Z(n)$.
The odd and even dihedral Gram formulas already hold at every
complex parameter. The four-bond fold is semisimple exactly when
$\delta(\delta-1)(\delta+2)\ne0$.

\subsection{Dependence on the polynomial form}
For example, consider the two free rank-two $\mathbb C[z]$-algebras
\[
 R=\mathbb C[z]\times\mathbb C[z],\qquad
 S=\mathbb C[z,x]/(x^2-z^2).
\]
After inverting $z$, evaluation at $x=z,-z$ gives an isomorphism
$S[z^{-1}]\to R[z^{-1}]$. Its inverse is
\[
 (a,b)\longmapsto a\frac{1+x/z}{2}+b\frac{1-x/z}{2}.
\]
At $z=0$, the fiber of $R$ is $\mathbb C\times\mathbb C$, which is semisimple.
The fiber of $S$ is $\mathbb C[x]/(x^2)$, whose nonzero radical $(x)$ has square zero.
We now construct an analogous pair of polynomial forms for the
triality-folded Brauer algebra of type $\mathscr{G}_2$.

\begin{proposition}\label{thm:g2-dual-zero-forms}
Let $R=\mathbb C[z]$, with loop parameter $z$, and
$L=\operatorname{Br}(\mathscr{G}_2;R[z^{-1}])$ in the normalization of
Proposition~\ref{prop:G2-criterion}. There are free rank-$39$
$R$-subalgebras $\mathcal G$ and $\widetilde{\mathcal G}$ of $L$
whose localizations are $L$, such that
\[
 \mathcal G(\delta)\text{ is semisimple}
 \Longleftrightarrow\delta\notin\{0,-2,1\},
\]
\[
 \widetilde{\mathcal G}(\delta)\text{ is semisimple}
 \Longleftrightarrow\delta\notin\{-2,1\}.
\]
Here the rescaled contraction is $\widetilde e_2=z^{-2}e_2$;
its square is $\widetilde e_2^{\,2}=z\widetilde e_2$.
\end{proposition}
\begin{proof}
Write $W=W(\mathscr{G}_2)$, $c=(r_1r_2)^3$, and choose
\[
 (u_1,u_2,u_3)=(1,r_2,r_1r_2),\qquad
 (v_1,v_2,v_3)=(1,r_1,r_2r_1).
\]
The normal-form basis of $L$ consists of the twelve elements of $W$,
the eighteen elements
$u_a e_1c^\epsilon u_b^{-1}$, and the nine elements
$m_{ab}=v_a e_2v_b^{-1}$, with $1\le a,b\le3$ and
$\epsilon\in\{0,1\}$. Let $S_1$ be the $R$-span of the middle
eighteen elements and $M$ the $R$-span of the $m_{ab}$. Define
\[
 \mathcal G=RW\oplus S_1\oplus M,\qquad
 \widetilde{\mathcal G}=RW\oplus S_1\oplus z^{-2}M
 \quad\text{inside }L.
\]
By the normal-form basis theorem \cite[Theorem~6.1]{LiuG22015}, both $\mathcal G$ and $\widetilde{\mathcal G}$ are free $R$-modules of rank $39$.

We now verify that both $\mathcal G$ and $\widetilde{\mathcal G}$ are closed under multiplication over $R$. The orbit of the central $\mathscr{D}_4$ simple root $\alpha_2$ under the
folded Weyl group $W=W(\mathscr{G}_2)$ consists of the three root lines
\[
 \alpha_2,\qquad
 \alpha_1+\alpha_2+\alpha_3+\alpha_4,\qquad
 \alpha_1+2\alpha_2+\alpha_3+\alpha_4,
\]
which are pairwise nonorthogonal. Using the relations
$e_1^2=ze_1$, $e_1r_2e_1=ce_1$, and $ce_1=e_1c$, we see that $S_1$ is closed under multiplication and no term from $M$ appears in $S_1\cdot S_1$. Moreover, $W$ permutes
the normal-form basis elements of $S_1$. By the frame action in the proof of
Proposition~\ref{prop:G2-criterion}, or directly from the relations
\[
 e_1e_2=r_2r_1e_2,\qquad e_2e_1=e_2r_1r_2,
\]
we see that left and right multiplication of $M$ by $RW\oplus S_1$ has
coefficients in $R=\mathbb C[z]$, and multiplying $M$ by $z^{-2}$ preserves these action coefficients. Finally, on $M$ we have
\[
 m_{ab}m_{cd}=z^2G(z)_{bc}m_{ad},\qquad
 G(z)=\begin{pmatrix}z&1&1\\1&z&1\\1&1&z\end{pmatrix}.
\]
Setting $\widetilde m_{ab}=z^{-2}m_{ab}$, we get
$\widetilde m_{ab}\widetilde m_{cd}=G(z)_{bc}\widetilde m_{ad}$.
Therefore, both $\mathcal G$ and $\widetilde{\mathcal G}$ are $R$-subalgebras of $L$ whose localizations at $z$ are equal to $L$.

At $z=0$, the two middle-layer cell forms are congruent to $G(0)$, which has
determinant $2\ne0$. Thus the quotient of both $\mathcal G(0)$ and $\widetilde{\mathcal G}(0)$ by the nine-dimensional minimal ideal is
semisimple. In $\mathcal G(0)$, the minimal ideal $M(0)$ satisfies $M(0)^2=0$, so $\mathcal G(0)$ is not semisimple. In $\widetilde{\mathcal G}(0)$, the minimal ideal is a $3\times3$ sandwich
matrix algebra with invertible pairing matrix $G(0)$, so it is isomorphic to $M_3(\mathbb C)$ and splits off as a direct factor by Lemma~\ref{lem:sandwich-layers}. Thus $\widetilde{\mathcal G}(0)$ is semisimple.
For $\delta\ne0$, both $\mathcal G(\delta)$ and $\widetilde{\mathcal G}(\delta)$ are isomorphic
to $\operatorname{Br}(\mathscr{G}_2;\delta)$, so Proposition~\ref{prop:G2-criterion} gives the result.
\end{proof}

For the remaining Laurent families, we make no assertion at zero
without a specified polynomial form. A relation such as $e_\kappa^2=0$
alone does not imply nonsemisimplicity: a full matrix algebra can contain
nonzero square-zero elements. The arguments above use nilpotent two-sided
ideals or the complete set of cell forms.

\section*{Computational data and reproducibility}
The accompanying verification package contains the programs, pinned
character-data files, and expected results for the finite calculations
in Section~\ref{sec:exceptional-cases}. The programs reconstruct the root
actions, geometric subgroup fusions, and layer matrices from their
specified inputs, and check the results by exact arithmetic. We use
SageMath~10.1; the package gives the execution commands and software versions,
and identifies the proof steps verified by each calculation.

The character values are obtained by specializing the included CHEVIE
sources at $q=1$ \cite{CHEVIE}. These character tables are fixed external
inputs; we compute the subgroup fusions from the geometric embeddings.
The checks use exact rational ranks and jet nullities for $\mathscr{E}_6$,
full-rank certificates over $\mathbb F_{1009}$ for the three residual
$\mathscr{E}_7$ matrices at $\delta=-7$, and exact annihilator and projector
identities for the remaining $\mathscr{E}_8$ channels.
For $\mathscr{F}_4,\mathscr{G}_2,\mathscr{H}_3,\mathscr{H}_4$, the package verifies the defining relations and the
determinant or spectral identities used in the proofs. The certification
principles are stated in Section~\ref{subsec:exact-certificate-principles}.
We include SHA-256 checksums and the notices for third-party data.
\ifdefempty{\verificationurl}{}{The verification code and data for this paper are available at \url{https://github.com/luanyongzhi/BrauerSemisimplicity}.}

\section*{Statements and Declarations}

\subsection*{Funding}
Y. Luan is supported by the Shandong Provincial Natural Science Foundation (Grant No. ZR2025MS36).

\subsection*{Competing interests}
The authors declare that they have no conflicts of interest.


\subsection*{Declaration of generative AI use}
During the preparation of this work, the authors used ChatGPT 6 Astra and Gemini 4 Argon for language editing, checking mathematical arguments and references, and simplifying proofs. The authors reviewed and edited the resulting content and take full responsibility for the manuscript.

\end{document}